\documentclass[12pt]{amsart}
\usepackage{amsmath}
\usepackage{amsxtra}
\usepackage{amscd}
\usepackage{amsthm}
\usepackage{amsfonts}
\usepackage{amssymb}
\usepackage{eucal}
\newcommand{\slth}{\widehat{\mathfrak{sl}}_2}

\newcommand{\wt}{{\rm wt}\,} 
 
\newcommand{\ch}{{\rm ch}\,}

\newcommand{\Vt}{V'}

\newcommand{\vb}{\overline{v}}
\newcommand{\sgn}{\mathop{\rm sgn}}
\newcommand{\Skew}{\mathop{\rm Skew}}
\newcommand{\seteq}{\mathbin{:=}}
\newcommand{\nn}{\nonumber}
\newcommand{\bea}{\begin{eqnarray}}
\newcommand{\ena}{\end{eqnarray}}
\newcommand{\be}{\begin{eqnarray*}}
\newcommand{\en}{\end{eqnarray*}}
\newcommand{\ba}{\begin{array}}
\newcommand{\ea}{\end{array}}
\newcommand{\C}{{\mathbb C}}
\newcommand{\Z}{{\mathbb Z}} 
\newcommand{\Q}{{\mathbb Q}} 
\newcommand{\I}{{\sqrt{-1}}}
\newcommand{\SG}{{\mathfrak{S}}}

\newcommand{\Cc}{\mathcal{C}}

\newcommand{\Zc}{\mathcal{Z}}

\newcommand{\Sym}{\mathop{\rm Sym}}

\newcommand{\F}{\mathcal{F}}

\newcommand{\Pc}{\mathcal{P}}
\newcommand{\Psit}{\widetilde{\Psi}}

\newcommand{\Ztc}{\widehat{\mathcal{Z}}}
\newcommand{\im}{{\sqrt{-1}}}
\newenvironment{tenumerate}{
  \begin{enumerate}
  
  }{\end{enumerate}}
\newcommand{\bi}{\begin{tenumerate}}
\newcommand{\ei}{\end{tenumerate}}
\newcommand{\isoto}[1][]%
{{\mathop{\buildrel{\sim}\over\longrightarrow}\limits_{#1}}}
\newcommand{\tpsi}{\tilde \psi}
\newcommand{\To}[1][\phantom{aaaa}]{\xrightarrow{\,#1\,}}

\numberwithin{equation}{section}
\newtheorem{thm}{Theorem}[section]
\newtheorem{prop}[thm]{Proposition}
\newtheorem{lem}[thm]{Lemma}
\newtheorem{cor}[thm]{Corollary}

\newcommand{\mono}{\rightarrowtail}
\newcommand{\ftens}{{\widehat{\mathop\otimes\limits_F}}}
\newcommand{\epi}{\twoheadrightarrow}
\begin{document} 
\title[A functional model for representations of
$U_q(\slth)$]
{A functional model for the tensor product
of level $1$ highest and level $-1$ lowest modules for
the quantum affine algebra $U_q(\slth)$}
\author{B. Feigin, M. Jimbo, M. Kashiwara, 
T. Miwa, E. Mukhin and Y. Takeyama}
\address{BF: Landau institute for Theoretical Physics, Chernogolovka,
142432, Russia}\email{feigin@feigin.mccme.ru}  
\address{MJ: Graduate School of Mathematical Sciences, The
University of Tokyo, Tokyo 153-8914, Japan}\email{jimbomic@ms.u-tokyo.ac.jp}
\address{MK: 
Research Institute for Mathematical Science, 
Kyoto University, Kyoto 606-8502,
Japan}\email{masaki@kurims.kyoto-u.ac.jp}
\address{TM: Department of Mathematics, Graduate School of Science,
Kyoto University, Kyoto 606-8502
Japan}\email{tetsuji@math.kyoto-u.ac.jp}
\address{EM: Department of Mathematics,
Indiana University-Purdue University-Indianapolis,
402 N.Blackford St., LD 270,
Indianapolis, IN 46202}\email{mukhin@math.iupui.edu}
\address{YT: Department of Mathematics, Graduate School of Science,
Kyoto University, Kyoto 606-8502
Japan}\email{takeyama@math.kyoto-u.ac.jp}

\date{\today}
\dedicatory{Dedicated to Alain Lascoux on the occasion 
of his sixtieth birthday}
\begin{abstract}
Let $V(\Lambda_i)$ (resp., $V(-\Lambda_j)$) be a fundamental integrable highest
(resp., lowest) weight module of $U_q(\slth)$. The tensor product
$V(\Lambda_i)\otimes V(-\Lambda_j)$ is filtered by submodules
$F_n=U_q(\slth)(v_i\otimes \vb_{n-i})$, $n\ge 0, n\equiv i-j\bmod 2$,
where $v_i\in V(\Lambda_i)$ is the highest vector and
$\vb_{n-i}\in V(-\Lambda_j)$ is an extremal vector. We show that
$F_n/F_{n+2}$ is isomorphic to the level $0$ extremal weight module 
$V(n(\Lambda_1-\Lambda_0))$. Using this we give a functional realization of 
the completion of $V(\Lambda_i)\otimes V(-\Lambda_j)$ by the filtration
$(F_n)_{n\geq0}$. The subspace of
$V(\Lambda_i)\otimes V(-\Lambda_j)$ of $\mathfrak{sl}_2$-weight $m$ 
is mapped to a certain space of sequences
$(P_{n,l})_{n\ge 0, n\equiv i-j\bmod 2,n-2l=m}$, 
whose members $P_{n,l}=P_{n,l}(X_1,\dots,X_l|z_1,\dots,z_n)$ are 
symmetric polynomials in $X_a$ and symmetric Laurent polynomials in $z_k$,
with additional constraints. When the parameter $q$ is specialized to $\I$, 
this construction settles a conjecture which arose in the study of form factors
in integrable field theory. 
\end{abstract}
\maketitle
\setcounter{section}{0}
\setcounter{equation}{0}
\section{Introduction}
For each fixed integer $m$ and $i\in\{0,1\}$,  
let us consider sequences 
${\bf p}=(P_{n,l})_{\substack{n\ge 0\\ n-2l=m}}$ 
of functions $P_{n,l}=P_{n,l}(X_1,\ldots,X_l|z_1,\ldots,z_n)$
satisfying the following conditions for all $n,l$: 
\begin{enumerate}
\item $P_{n,l}$ is a
polynomial in $X_1,\ldots,X_l$ which is 
skew-symmetric when $l>1$, 
\item $P_{n,l}$ is a symmetric Laurent polynomial in 
$z_1,\ldots,z_n$,
\item $\deg_{X_a}P_{n,l}\le n-1$,
\item 
\be
&&P_{n+2,l+1}(X_1,\ldots,X_l,z^{-1}|z_1,\ldots,z_n,z,-z)
\\
&&\quad
=z^{-n-1+i}\prod_{a=1}^l(1-X_a^2z^2)
\cdot P_{n,l}(X_1,\ldots,X_l|z_1,\ldots,z_n).
\en
\end{enumerate}
Such sequences naturally arise in the form factor 
bootstrap approach 
to massive integrable models of quantum field theory \cite{Sbk}.
Form factors are sequences of 
matrix elements of local fields taken 
between the vacuum and the asymptotic states. 
They are typically given by certain integrals 
involving polynomials $P_{n,l}$ of the type mentioned above, 
wherein $\alpha_a=-\log X_a$ are the integration variables 
and $\beta_j=\log z_j$ are the rapidity variables of 
asymptotic particles. 
More specifically, 
the sequences ${\bf p}$ satisfying (i)--(iv) 
appear in the sine-Gordon and the 
$SU(2)$ invariant Thirring models, 
and are called `$\infty$-cycles' of weight $m$ in \cite{JMMT}. 
The conditions (i)--(iii), along with (iv), which was originally
proposed in \cite{NT} (see \cite{JMMT} in the present form with $i=0$),
ensure the locality of the fields.

Denote by $\Ztc_{\C}^{{\rm skew}(i,j)}[m]$ the space of
all $\infty$-cycles of weight $m$ with $m\equiv i-j\bmod 2$. 
It was shown in \cite{JMMT}
\footnote{In \cite{JMMT}, 
the space $\Ztc_{\C}^{{\rm skew}(0,j)}$ was
denoted by $\Ztc^{(j)}$.}
 that the space
$\Ztc_{\C}^{{\rm skew}(i,j)}:=\oplus_{\substack{m\in\Z\\m\equiv i-j}}
\Ztc_{\C}^{{\rm skew}(i,j)}[m]$ 
admits an action of the quantum affine algebra 
$U_{\im}(\slth)$ with the parameter $q=\im$. 
It was conjectured further that $\Ztc_{\C}^{{\rm skew}(i,j)}$
is isomorphic to the tensor product module 
$V_{\im}(\Lambda_i)\otimes V_{\im}(-\Lambda_j)$ 
(with a proper completion, see below) 
of integrable modules of level $1$ and level $-1$, respectively.  
The purpose of this paper is to clarify the representation 
theoretical origin of $\infty$-cycles, 
and to supply a proof of the above conjecture. 

Though only the case $q=\im$ is relevant to form factors, analogs of
$\infty$-cycles exist also for generic $q$. 
In the below we outline their construction. 
Let $U_q=U_q(\slth)$ be the quantum affine algebra over 
$K=\C(q)$, and let $U'_q$ be the subalgebra with the scaling element being 
dropped. For $i,j=0,1$, let $V(\Lambda_i)$ (resp., $V(-\Lambda_j)$)
be the integrable highest weight (resp., lowest weight) module
with highest weight $\Lambda_i$ (resp., $-\Lambda_j$) and
highest weight vector $v_i$ (resp., lowest weight vector $\bar{v}_{-j}$).
Let further $\bar{v}_n\in V(-\Lambda_j)$ 
($n\equiv j\bmod 2$) be an extremal vector obtained from 
$\bar{v}_{-j}\in V(-\Lambda_j)$ by the braid group action
corresponding to the translation element  
$(s_0s_1)^{-(n+j)/2}$ of the Weyl group. 
In the tensor product $V(\Lambda_i)\otimes V(-\Lambda_j)$, the submodules 
\be
F^{(i,j)}_n=U_q(v_i\otimes \bar{v}_{n-i})
\qquad
(n\ge 0, n\equiv i-j\bmod 2)
\en
define a decreasing filtration 
\bea
V(\Lambda_i)\otimes V(-\Lambda_j)
=F^{(i,j)}_{|i-j|}\supset
\cdots\supset F^{(i,j)}_n\supset F^{(i,j)}_{n+2}\supset\cdots.
\label{eqn:filter}
\ena
Denote by $V_z=V\otimes K[z,z^{-1}]$ the evaluation module 
based on the two-dimensional space $V=K v_+\oplus K v_-$. 
Then there exists a $U'_q$-linear map
\be
\psi_n~:~V(\Lambda_i)\otimes V(-\Lambda_j)
\longrightarrow(V_{z_1}\otimes\cdots\otimes V_{z_n})^\wedge
\en
such that $\psi_n(v_i\otimes \bar{v}_{n-i})=v_+^{\otimes n}$
and  
$\psi_n(F^{(i,j)}_{n+2})=0$.
Here the right hand side means the completion 
$\bigl(V_{z_1}\otimes\cdots\otimes V_{z_n}\bigr)
\otimes_{K[z_1/z_2,\ldots,z_{n-1}/z_n]} K[[z_1/z_2,\ldots,z_{n-1}/z_n]]$. 
Furthermore, $\psi_n$ induces an isomorphism 
\bea
\phi_n~:~F^{(i,j)}_n/F^{(i,j)}_{n+2}
\overset{\sim}{\longrightarrow}
V(n(\Lambda_1-\Lambda_0))
\label{eqn:FFV}
\ena
between the associated graded space 
$F^{(i,j)}_n/F^{(i,j)}_{n+2}$ and
the extremal weight module 
$V(n(\Lambda_1-\Lambda_0))=U_qv_+^{\otimes n}$ of level $0$
(Theorem \ref{thm:FFV}). 

For $0\le l\le n$, let $\F_{n,l}$ denote the space of 
symmetric polynomials $P(X_1,\ldots,X_l)$, 
with coefficients in $K[z_1^{\pm 1},\ldots,z_n^{\pm 1}]$, 
such that $\deg_{X_a}P\le n-1$ for each $a$ and 
\be
P|_{X_1=q^{-2}X_2=z_j^{-1}}=0\hbox{ for each }j=1,\ldots,n \quad \mbox{when $l>1$.}
\en
Let $\F_n=\oplus_{l=0}^n\F_{n,l}$. 
It is known \cite{TV} that there is an embedding 
of the tensor product of evaluation modules 
\be
\Cc_n~:~V_{z_1}\otimes\cdots\otimes V_{z_n}\longrightarrow 
\F_n,
\en 
which is $K[z_1,\ldots,z_n]$-linear.
The subspace $(V_{z_1}\otimes\cdots\otimes V_{z_n})_m$
of $\mathfrak{sl}_2$-weight $m$ is mapped to $\F_{n,l}$ where $n-2l=m$.
Let $\widehat{\Cc}_n:(V_{z_1}\otimes\cdots\otimes V_{z_n})^\wedge\rightarrow
\widehat\F_n$ be the extension of $\Cc_n$ where
$\widehat\F_n=\F_n\otimes_{K[z_1/z_2,\ldots,z_{n-1}/z_n]} K[[z_1/z_2,\ldots,z_{n-1}/z_n]]$.
This is an isomorphism. 
Analogous isomorphisms hold also for modules over the integral form 
$U_A\subset U_q$ (where $A=\C[q,q^{-1}]$).
It turns out that the image of 
$\varphi_n=\widehat{\Cc}_n\circ\psi_n$  
is contained in the $\mathfrak{S}_n$-invariant subspace $\F_n^{\mathfrak{S}_n}$
of $\F_n$ (without completion). Here the symmetric group $\mathfrak{S}_n$ acts on $\F_n$
by the permutation of the variables $z_1\ldots,z_n$.
The image of the map 
\be
\varphi=\prod_{n}
\varphi_{n}
~:~V(\Lambda_i)\otimes V(-\Lambda_j)
\longrightarrow \prod_{\substack{n\ge 0\\n\equiv i-j\bmod 2}}\F_n^{\mathfrak{S}_n}
\en
is contained in the subspace defined as follows.
Denote by $\Ztc^{(i,j)}[m]$ the space of all sequences 
${\bf p}=(P_{n,l})_{\substack{n\ge 0\\ n-2l=m}}$ 
of polynomials $P_{n,l}\in\F_{n,l}^{\mathfrak{S}_n}$,  
satisfying the property:
\be
{\rm (iv)'}&&
\mbox{
$P_{n+2,l+1}(X_1,\ldots,X_l,z^{-1}|z_1,\ldots,z_n,z,q^2z)$
}
\\
&&\quad 
\mbox{$=
z^{-n-1+i}\prod_{a=1}^l(1-q^{-2}X_az)(1-q^{2}X_az)\cdot
P_{n,l}(X_1,\ldots,X_l|z_1,\ldots,z_n)$}.
\en
We set $\Ztc^{(i,j)}=\oplus_{m\in\Z}\Ztc^{(i,j)}[m]\subset
\prod_{\substack{n\ge 0\\n\equiv i-j\bmod 2}}\F_n^{\mathfrak{S}_n}$.
The image of $V(\Lambda_i)\otimes V(-\Lambda_j)$ is contained in this subspace,
and moreover, the completion of $V(\Lambda_i)\otimes V(-\Lambda_j)$
by the filtration $\{F^{(i,j)}_n\}$ is isomorphic to $\Ztc^{(i,j)}$
(Theorem \ref{thm:functional-realization}). 

By a simple transformation, the specialization of 
(the integral form of) $\Ztc^{(i,j)}$ to $q=\im$  
is mapped injectively to the space of 
$\infty$-cycles $\Ztc_{\C}^{{\rm skew}(i,j)}$. From this follows the
conjectured isomorphism in the original setting.

Quite generally, it is known \cite{BN} 
for an arbitrary quantized affine algebra that 
the tensor product of highest and lowest modules with 
total level zero  
admits a filtration with a property 
similar to \eqref{eqn:FFV}.  
However each filter is in general not generated by tensor products of 
extremal vectors.
It would be interesting to study their structure. 
In particular the filtration on the tensor product induces 
a filtration on the lowest weight module. 
In the case of $\slth$ we give a conjecture 
on the character of the associated graded space for the latter 
(see \eqref{MEL}). 


The text is organized as follows. 
In Section 2, we give a brief review on extremal weight
modules and set up the notation. 
We then 
introduce the filtration \eqref{eqn:filter}
of the tensor product of level $1$ highest and level $-1$ lowest modules, 
and prove the isomorphism \eqref{eqn:FFV}. 
In Section 3, we discuss the polynomial realization of the 
tensor product $V(\Lambda_i)\otimes V(-\Lambda_j)$ 
and the associated graded spaces of the filtration $\{F^{(i,j)}_n\}$. 
The main results are stated in 
Theorem \ref{thm:functional-realization} and 
Theorem \ref{thm:VVZ}. 

For the reader's convenience, 
we summarize in Appendix 
some basic facts concerning crystal and 
global basis of extremal weight modules used in the text. 
We also give a brief account of the filtration 
of the tensor product modules $V(\xi)\otimes V(-\eta)$ 
for general quantum affine algebras. 

\section{Filtration by extremal vectors}\label{sec:filt}
\subsection{Notation}\label{subsec:notation}
First we fix our notation concerning quantum affine algebra $U_q(\slth)$. 
Set $I=\{0,1\}$. 
Let $P=\Z \Lambda_0\oplus\Z \Lambda_1\oplus\Z\delta$ be 
the weight lattice for $\widehat{\mathfrak{sl}}_2$,
$P^*=\Z h_0\oplus\Z h_1\oplus\Z d$ its dual lattice and
$P_+=\{\lambda\in P|\langle\lambda,h_i\rangle\geq0\,
(i\in I)\}$.
The quantum affine algebra 
$U_q=U_q(\widehat{\mathfrak{sl}}_2)$ is the 
algebra over $K=\C(q)$ generated by $e_i,f_i$ 
$(i\in I)$ and $q^h$ ($h\in P^*$),  
under the defining relations 
\be
&&q^hq^{h'}=q^{h+h'},\quad q^0=1,
\\
&&q^he_iq^{-h}=q^{\langle{h,\alpha_i\rangle}}e_i,
\quad q^hf_iq^{-h}=q^{-\langle{h,\alpha_i\rangle}}f_i,
\\
&&[e_i,f_j]=\delta_{ij}\frac{t_i-t_i^{-1}}{q-q^{-1}},
\\
&&\sum_{r=0}^3(-1)^re_i^{(3-r)}e_je_i^{(r)}=0
\quad (i\neq j),\\
&&\sum_{r=0}^3(-1)^rf_i^{(3-r)}f_jf_i^{(r)}=0
\quad (i\neq j).\\
\en
Here $t_i=q^{h_i}$, $\alpha_1=2(\Lambda_1-\Lambda_0)$, 
$\alpha_0=\delta-\alpha_1$, and for an element $x\in U_q$,
we denote by $x^{(r)}$ the divided power $x^r/[r]!$, 
where $[r]!=\prod_{j=1}^r[j]$, $[j]=(q^j-q^{-j})/(q-q^{-1})$.
The element $C\seteq t_0t_1$ is central, and 
$D\seteq q^d$ is the scaling element. 
We will use the coproduct
\bea
\Delta(e_i)=e_i\otimes t_i^{-1}+1\otimes e_i,~
\Delta(f_i)=f_i\otimes 1+t_i\otimes f_i,~
\Delta(q^h)=q^h\otimes q^h. 
\label{copro}
\ena

We denote by $W$ the Weyl group for $\widehat{\mathfrak{sl}}_2$.

We say that a $U_q$-module $M$ is integrable if the action of 
$e_i,f_i$ ($i\in I$) is locally nilpotent 
and $M=\oplus_{\mu\in P}M_\mu$, 
$M_\mu\seteq\{u\in M \mid q^hu=q^{\langle h,\mu\rangle} u~~(h\in P^*)\}$.
For $u\in M_\mu$ we write $\wt u=\mu$, $D u=q^{\deg u}u$.

We will consider the specialization at $q=\sqrt{-1}$.
For this purpose, we need the integral form of $U_q$.
We set $A=\C[q,q^{-1}]$. Let $U_A$ be the $A$-subalgebra of $U_q$ generated by 
$e_i^{(r)},f_i^{(r)}$ $(i\in I,\,r\in\Z_{\geq0})$ and $q^h$ $(h\in P^*)$. 
For $\epsilon\in\C\backslash\{0\}$, the specialization $U_\epsilon$ is 
the $\C$-algebra $U_A/U_A(q-\epsilon)$.
We denote by $U_q^+$ (resp., $U_q^-$) the subalgebra of $U_q$ generated by
$e_i$ $(i\in I)$ (resp., $f_i$ $(i\in I)$).

We denote by $U_q^{\geq0}$ the subalgebra of $U_q$ generated by 
$e_i$, $q^h$ ($i\in I,h\in P^*$) and $f_1$, and by
$U^{\geq0}_A$ the $A$-subalgebra generated by $e^{(r)}_i$, $q^h$ and 
$f^{(r)}_1$ for $i\in I,h\in P^*$, $r\in\Z_{\ge 0}$. 
Likewise, we define $U_q^{\leq0}$ and $U_A^{\leq0}$ by 
changing $e_i$ (or $e_i^{(r)}$) to $f_i$ (or $f_i^{(r)}$) and $f_1$
(or $f_1^{(r)}$) to $e_1$ (or $e_1^{(r)}$)
in the definition above.

Let $M$ be a $K$-vector space.
An $A$-submodule $M_A$ of $M$ is called an $A$-lattice of $M$ if it is a free
$A$-module and $M=M_A\otimes_AK$.
For an $A$-lattice $M_A$ and $\epsilon\in\C\backslash\{0\}$,
we write $(M_A)_\epsilon=M_A/(q-\epsilon)M_A$.
We call it the specialization of $M_A$ at $q=\epsilon$. If
there is no fear of confusion, we abbreviate $(M_A)_\epsilon$ to $M_\epsilon$.
When we specialize a $U_q$-module $M$ to $q=\epsilon$, we must first choose
an $A$-lattice of $M$ which is 
stable under the action of $U_A$.
Then, the specialization $M_\epsilon$ admits a $U_\epsilon$-action and
we obtain a $U_\epsilon$-module. We note that different $A$-lattices
may lead to non-isomorphic $U_\epsilon$-modules.

Let $M$ be an integrable $U_q$-module.
Consider the bi-grading of $M=\oplus_{a,n\in\Z}M_{a,n}$,
where $M_{a,n}=\{u\in M\mid Du=q^au, t_1u=q^nu\}$.
We define its character by 
\bea\label{CHARDEF}
\ch_{v,z}M\seteq \sum_{a,n\in\Z}\dim_KM_{a,n}v^az^n.
\ena
If $M$ is a bi-graded $\C$-vector space (resp., a bi-graded free $A$-module)
we define its character by (\ref{CHARDEF})
replacing the dimension over $K$ by the dimension over $\C$
(resp., the rank over $A$).
If $M_A$ is an $A$-lattice of $M$, the characters 
$\ch_{v,z}M$, $\ch_{v,z}M_A$, $\ch_{v,z}M_\epsilon$ are all equal.
We note that the character is well-defined only if each subspace
$M_{a,n}$ is finite-dimensional (or of finite rank). In fact,
$U_q$-modules we consider in this paper do not necessarily satisfy
this property, e.g., $V(2(\Lambda_1-\Lambda_0))$ given in the next section.
\subsection{Extremal weight modules}\label{subsec:etremal-mod}
We recall the notion of extremal weight 
modules over $U_q$, introduced in \cite{K1}
for general quantized enveloping algebras. 
Let $\lambda\in P$. In the present case of affine type algebras, 
the extremal weight module $V(\lambda)$ is characterized as
the universal integrable $U_q$-module with the following defining relations:
\bea
V(\lambda)=U_q u_\lambda\hbox{ where }{\rm wt}\,u_\lambda=\lambda,\\
{\rm wt}\,V(\lambda)\subset\hbox{the convex hull of }W\lambda.
\ena
We define an $A$-lattice of $V(\lambda)$ by
$V_A(\lambda)=U_A u_\lambda$. We denote its specialization to
$q=\epsilon\in\C\backslash\{0\}$ by $V_\epsilon(\lambda)$.

The extremal weight module $V(\lambda)$ has 
a family of extremal vectors 
$S_wu_\lambda$ indexed by $w\in W$.  
They are defined by $S_{\rm id}u_\lambda=u_\lambda$ and
\bea
S_{s_iw}u_\lambda=
\begin{cases}
f_i^{(\langle h_i,w\lambda\rangle)}S_wu_\lambda
&\hbox{ if $\langle h_i,w\lambda\rangle\geq0$,}\\[5pt]
e_i^{(-\langle h_i,w\lambda\rangle)}S_wu_\lambda
&\hbox{ if $\langle h_i,w\lambda\rangle\leq0$.}
\end{cases}
\ena
For $w\in W$, we have $V(\lambda)=U_qS_wu_\lambda$. 
There is a canonical isomorphism
\bea\label{ISOMOR}
V(w\lambda)=U_qu_{w\lambda}\isoto V(\lambda)
\ena
sending $u_{w\lambda}$ to $S_wu_\lambda$.

If $\lambda\in P_+$, $V(\lambda)=U_qu_\lambda$ is the 
integrable highest weight module with highest weight $\lambda$ and
highest weight vector $u_\lambda$. Similarly,
$V(-\lambda)=U_qu_{-\lambda}$ is the integrable lowest weight module
with lowest weight $-\lambda$ and lowest weight vector $u_{-\lambda}$.
The following result \cite{L} is basic in our study.
\begin{prop}\label{Lu}
Let $\lambda,\mu\in P_+$ and $\epsilon\in\C\backslash\{0\}$.  
The tensor products $V(\lambda)\otimes V(-\mu)$,
$V_A(\lambda)\otimes V_A(-\mu)$, $V_\epsilon(\lambda)\otimes V_\epsilon(-\mu)$
have the cyclic vector $u_\lambda\otimes u_{-\mu}$.
They are characterized as the universal cyclic module 
with the cyclic vector $v$,  
with weight condition 
$\wt v=\lambda-\mu$ and the defining relations
\bea
f_i^{(r)}v=0\hbox{ for any }r\geq\langle h_i,\lambda\rangle+1,\\
e_i^{(r)}v=0\hbox{ for any }r\geq\langle h_i,\mu\rangle+1.
\ena
\end{prop}
In this paper, we specifically consider the $U_q$-modules
$V(\Lambda_i)$ and $V(-\Lambda_i)$ where $i=0,1$. Let
$v_i\seteq u_{\Lambda_i}\in V(\Lambda_i)$ be the highest weight vector,
and $\vb_{-i}\seteq u_{-\Lambda_i}\in V(-\Lambda_i)$ the lowest weight vector. 
For $n\in\Z$ such that $n\equiv i\bmod 2$, consider the extremal vectors
\bea
&&v_n\seteq S_{(s_0s_1)^{(n-i)/2}}v_i,
\quad \wt v_n=\Lambda_0+n(\Lambda_1-\Lambda_0)-\frac{n^2-i}{4}\delta,
\label{EXTDEG1}\\
&&\vb_n\seteq S_{(s_0s_1)^{(-n-i)/2}}\vb_{-i},
\quad \wt \vb_n=-\Lambda_0+n(\Lambda_1-\Lambda_0)+
\frac{n^2-i}{4}\delta.\label{EXTDEG2}
\ena

We shall deal also with level $0$ weights
$\lambda=n(\Lambda_1-\Lambda_0)+r\delta$ ($n,r\in\Z$). 
Let $U'_q$ be the $K$-subalgebra of $U_q$ generated by
$e_i,f_i,t_i^{\pm1}$ $(i\in I)$. 
We have an isomorphism
$V(\lambda+m\delta)\simeq V(\lambda)$ as $U'_q$-modules for any
$m\in\Z$. Indeed, as $U'_q$-modules an isomorphism is defined by sending
$u_{\lambda+m\delta}$ to $u_\lambda$. If $m\not=0$, this map is not $D$-linear.
The degree $d$ in $V(\lambda+m\delta)$ and that in $V(\lambda)$ differs by
$m$. From the structure of the modules, which we will describe
below, it is easy to see that there exists an isomorphism
as $U_q$-modules between
$V(\lambda+m\delta)$ and $V(\lambda)$ if and only if $m\equiv0\bmod n$.

Set $V=K v_+\oplus K v_-$, $V_z=V\otimes K[z^{\pm1}]$. 
We abbreviate $v_\pm\otimes z^m$ to $z^mv_\pm$
and set $D z^mv_\pm=q^m z^mv_\pm$. 
Define further an action of $U'_q$ on $V_z$ by the rule that it
commutes with the multiplication by $z$, and 
\bea
&&
\ba{lllll}
t_0v_\pm=q^{\mp1}v_\pm,&
e_0v_+=zv_-,& e_0v_-=0,&f_0v_+=0,& f_0v_-=z^{-1}v_+,\\
t_1v_\pm=q^{\pm1}v_\pm,&
e_1v_+=0,& e_1v_-=v_+,& f_1v_+=v_-,& f_1v_-=0.
\ea
\ena
Then we have an isomorphism of $U_q$-modules
$V_z\simeq V(\Lambda_1-\Lambda_0)$ 
which sends $v_+$ to $u_{\Lambda_1-\Lambda_0}$.  
We have also
$V_{A,z}\simeq V_A(\Lambda_1-\Lambda_0)$,  
where
$V_{A,z}=(Av_+\oplus Av_-)\otimes_A  A[z^{\pm1}]$. 
In addition, $S_{(s_0s_1)^m}v_+=z^{-m}v_+$ and 
$S_{s_1(s_0s_1)^m}v_+=z^{-m}v_-$.

For general $n\in\Z_{\ge 0}$, we identify  $V_z^{\otimes n}$
with the vector space $V^{\otimes n}\otimes\ K[z_1^{\pm1},\ldots,z_n^{\pm1}]$.
The $U_q$ module $V(n(\Lambda_1-\Lambda_0))$ is isomorphic to a submodule
$U_qv_+^{\otimes n}$ of $V_z^{\otimes n}$. We identify them also,
by the identification $v_+^{\otimes n}=u_{n(\Lambda_1-\Lambda_0)}$.
Denoting by $\SG_n$ the symmetric group on $n$ letters, we have 
\bea
\oplus_{m\in\Z}V(n(\Lambda_1-\Lambda_0))_{n(\Lambda_1-\Lambda_0)+m\delta}
=K v_+^{\otimes n}
\otimes K[z_1^{\pm 1},\ldots,z_n^{\pm1}]^{\SG_n}.\label{NISO}
\ena

Viewed as a subspace of 
$V^{\otimes n}\otimes\ K[z_1^{\pm1},\ldots,z_n^{\pm1}]$,
$V(n(\Lambda_1-\Lambda_0))$ is invariant under multiplication
by symmetric Laurent polynomials, and the multiplication commutes
with the $U'_q$-action. We have
$S_{(s_0s_1)^m}v_+^{\otimes n}=z^{-m}v_+^{\otimes n}$ and
$S_{s_1(s_0s_1)^m}v_+^{\otimes n}=z^{-m}v_-^{\otimes n}$
where $z=z_1\cdots z_n$.

The $U_q$-module $V(\lambda)$ with $\lambda=n(\Lambda_1-\Lambda_0)$
is also characterized as the universal integrable module
$U_qu_\lambda$ satisfying the following properties.
\bea
&&{\rm wt}\,u_\lambda=\lambda,
\label{level0-1}\\
&&V(\lambda)_\xi=0\hbox{ if }\xi\in \lambda+\Z_{>0}\alpha_1+\Z\delta.
\label{level0-2}
\ena
\subsection{Filtration of $V(\Lambda_i)$}
We introduce a decreasing filtration of $V(\Lambda_i)$ by
$U_q^{\leq0}$-modules. Namely, we set
\be
F^{(i)}_n\seteq U_q^{\leq0}v_n\quad (n\ge 0, n\equiv i\bmod 2).
\en
Then, we have $F^{(i)}_n\supset F^{(i)}_{n+2}$ and
$\cap_{n}F^{(i)}_n=0$.
Similarly, we define the filtration of $V(-\Lambda_i)$ by
$U_q^{\geq0}$-modules.
\be
\overline F^{(i)}_n\seteq U_q^{\geq0}\vb_n\quad (n\ge 0, n\equiv i\bmod 2).
\en
\subsection{Filtration of $V(\Lambda_i)\otimes V(-\Lambda_j)$}
\label{subsec:filt}
Hereafter we fix $i,j\in I$ and consider the tensor product 
$V(\Lambda_i)\otimes V(-\Lambda_j)$. 
For $n\ge 0$ such that $n\equiv i-j\bmod 2$, we set 
\be
F^{(i,j)}_n\seteq U_q(v_i\otimes \vb_{n-i}).
\en
Note that $F^{(i,j)}_n=U_q(v_{-n+j}\otimes\vb_{-j})$.
This defines a decreasing filtration by $U_q$-submodules 
\be
V(\Lambda_i)\otimes V(-\Lambda_j)=F^{(i,j)}_{|i-j|}
\supset F^{(i,j)}_{|i-j|+2}\supset\cdots. 
\en
Similarly, we set
$F^{(i,j)}_{A,n}\seteq U_A(v_i\otimes\vb_{n-i})$,
$F^{(i,j)}_{\epsilon,n}\seteq U_\epsilon(v_i\otimes\vb_{n-i})$.

\begin{prop}\label{prop:intersection}
\be
&&\bigcap_{\substack{n\ge 0\\n\equiv i-j \bmod 2}} F^{(i,j)}_n=0.
\en
\end{prop}
A proof is given in Appendix.

Now we state one of the main results in this paper.
\begin{thm}\label{thm:FFV}
We have an isomorphism of $U'_q$-modules
\be
&&\phi_n:F^{(i,j)}_n/F^{(i,j)}_{n+2}\simeq 
V(n(\Lambda_1-\Lambda_0))\quad(n\ge 0, n\equiv i-j\bmod 2),
\en
which sends $v_i\otimes\vb_{n-i}$ to $v_+^{\otimes n}$.
For a weight vector $u\in F^{(i,j)}_n/F^{(i,j)}_{n+2}$,
we have
\be
\deg u=\deg\phi_n(u)+\frac{(n-i)^2-j}4.
\en
\end{thm}
The proof is based on the following two propositions. 

\begin{prop}\label{prop:support}
We have 
\be
(F^{(i,j)}_n/F^{(i,j)}_{n+2})_\xi=0\quad
\mbox{{\it for} }\quad
\xi\in n(\Lambda_1-\Lambda_0)+\Z_{>0}\alpha_1+\Z\delta. 
\en
\end{prop}

\begin{prop}\label{prop:FtoV}
There exists a $U'_q$-linear surjection
\be
\tpsi_n:F^{(i,j)}_n\longrightarrow V(n(\Lambda_1-\Lambda_0)) 
\en 
which sends $v_i\otimes\vb_{n-i}$ to $v_+^{\otimes n}$.
\end{prop}

We show Proposition \ref{prop:support} in 
Sections \ref{subsec:level1}--\ref{subsec:level1-2}, 
and Proposition \ref{prop:FtoV} in Section \ref{subsec:intertwiners}. 
Assuming them, let us prove Theorem \ref{thm:FFV}.
\medskip

\noindent
{\it Proof of Theorem \ref{thm:FFV}.}\quad 
Set $\lambda=n(\Lambda_1-\Lambda_0)$. 
The module $F^{(i,j)}_n/F^{(i,j)}_{n+2}$ is generated by a vector
$u\seteq v_i\otimes \vb_{n-i}$ of weight $\lambda+((n-i)^2-j)/4\cdot\delta$, 
and the weight spaces are restricted
as in Proposition \ref{prop:support}.
{}From the characterization 
\eqref{level0-1}, \eqref{level0-2} 
of extremal weight modules $V(\lambda)$, 
we have a surjection of $U'_q$-modules
\bea
V(\lambda) \longrightarrow F^{(i,j)}_n/F^{(i,j)}_{n+2}. 
\label{map1}
\ena
Since $\wt\psi_n(v_i\otimes\vb_{n+2-i})+\Z\delta$ does not appear 
in the weights of $V(\lambda)$, we have $\tilde{\psi}_n(F^{(i,j)}_{n+2})=0$. 
Hence we have also a surjection 
\bea
F^{(i,j)}_n/F^{(i,j)}_{n+2}\longrightarrow V(\lambda).
\label{map2}
\ena
The mappings \eqref{map1} and \eqref{map2} 
exchange the cyclic vector $u_\lambda\in V(\lambda)$ and the cyclic vector
$u\in F^{(i,j)}_n/F^{(i,j)}_{n+2}$.
These maps are therefore $U'_q$-isomorphisms. The statement about the degree
follows from (\ref{EXTDEG1}) and (\ref{EXTDEG2}).
\qed
\begin{prop}\label{CORO}
The isomorphism $\phi_n$ in Theorem \ref{thm:FFV} induces the following isomorphisms{\/\rm :}
\be
F^{(i,j)}_{A,n}/F^{(i,j)}_{A,n+2}\simeq V_A(n(\Lambda_1-\Lambda_0)),\\
F^{(i,j)}_{\epsilon,n}/F^{(i,j)}_{\epsilon,n+2}
\simeq V_\epsilon(n(\Lambda_1-\Lambda_0)).
\en
\end{prop}
A proof is given in Appendix.

The following result is a consequence of Theorem \ref{thm:FFV}.
We use a result in \cite{K3} in the proof.
\begin{prop}\label{prop:BVN}
\bea\label{BVN}
\overline F^{(j)}_n/\overline F^{(j)}_{n+2}
\simeq U^{\geq0}_qv_+^{\otimes n}.
\ena
\end{prop}
\begin{proof}
we show that
\be
(v_0\otimes\overline F^{(j)}_n)\cap U_q(v_0\otimes\vb_{-(n+2)})
=v_0\otimes\overline F^{(j)}_{n+2}.
\en
The inclusion $\supset$ is clear. In order to show the other inclusion,
we apply Corollary 3.2 in \cite{K3}
by setting $\lambda=\Lambda_0$, $\mu={\rm wt}\,\vb_{-(n+2)}$.
Set $M=\overline F^{(j)}_n=U^+_q\vb_{-n}$ and
$N=\overline F^{(j)}_{n+2}=U^+_q\vb_{-(n+2)}$.
Note that $v_0=u_\lambda$ and $V(\mu)\simeq V(-\Lambda_j)$. We have
\be
(u_\lambda\otimes M)\cap U_q(u_\lambda\otimes\vb_{-(n+2)})
\subset(u_\lambda\otimes V(\mu))\cap U_q(u_\lambda\otimes\vb_{-(n+2)})
=u_\lambda\otimes N.
\en
Noting that
$U_q(v_0\otimes \vb_{n+2})=U_q(v_0\otimes \vb_{-(n+2)})$, we have
\be
M/N\simeq v_0\otimes M/((v_0\otimes M)\cap U_q(v_0\otimes \vb_{n+2}))
\subset F^{(0,j)}_n/F^{(0,j)}_{n+2}\simeq V(n(\Lambda_1-\Lambda_0)).
\en
Since $v_0\otimes\vb_n$ is identified with $v_+^{\otimes n}$
in $V(n(\Lambda_1-\Lambda_0))\subset V_z^{\otimes n}$, we have (\ref{BVN}).
\end{proof}
See Corollary \ref{CHAR+} for the character of \eqref{BVN}.
{}From the proof of Proposition \ref{prop:BVN} 
we have also 
\bea 
F_{n}^{(0, j)} \cap (v_{0} \otimes V(-\Lambda_{j}))=
v_{0} \otimes \overline{F}_{n}^{(j)}. 
\label{induced-filtration} 
\ena

\subsection{Vertex operator realization}\label{subsec:level1}
To show Proposition \ref{prop:support}, we utilize the
realization of the modules $V(\Lambda_i)$, $V(-\Lambda_j)$ 
in terms of vertex operators. 
For that purpose, consider the Drinfeld generators 
of $U_q$,   
$x^{\pm}_k$ ($k\in\Z$), $b_n$ ($n\in\Z\backslash\{0\}$) 
and $t_1^{\pm 1}$, $C^{\pm 1}$, $D^{\pm 1}$ 
obeying the following relations
\footnote{These generators are the same as those in \cite{JM}, except 
$C=\gamma$, $D=q^d$ and $b_n=(n/[n])\gamma^{n/2}a_n$.}. 
\bea
&&Db_nD^{-1}=q^{n}b_n,~~
Dx^{\pm}_nD^{-1}=q^{n}x^{\pm}_n,
\label{Dr1}\\
&&t_1b_n=b_nt_1,\quad 
t_1x^{\pm}_nt_1^{-1}=q^{\pm 2}x^{\pm}_n,
\label{Dr2}
\\
&&
[b_m,b_n]=m\frac{[2m]}{[m]^2}
\frac{C^{m}-C^{-m}}{q-q^{-1}}\delta_{m+n,0},
\label{Dr3}
\\
&&
[b_n,x_k^{\pm}]=\pm
\frac{[2n]}{[n]}C^{(n\mp|n|)/2}x^{\pm}_{k+n},
\label{Dr4}
\\
&&x^{\pm}_{k+1}x^{\pm}_l-q^{\pm 2}x_l^{\pm}x_{k+1}^\pm
=q^{\pm 2}x^{\pm}_{k}x^{\pm}_{l+1}-x_{l+1}^{\pm}x_{k}^\pm,
\label{Dr5}
\\
&&[x^+_k,x_l^-]=
\frac{C^{-l}\varphi^+_{k+l}-C^{-k}\varphi^{-}_{k+l}}{q-q^{-1}}.
\label{Dr6}
\ena Here \be \sum_{k\in \Z}\varphi^{\pm}_{\pm k}z^k =t_1^{\pm
1}\exp\left(\pm\sum_{n=1}^\infty\frac{(q^n-q^{-n})}{n} b_{\pm
n}z^n\right).  
\en 
They are related to the  Chevalley generators by 
\footnote{
To conform with the coproduct \eqref{copro}, we have changed
the identification slightly from \cite{JM}. With this identification,
the formulas of coproduct for the Drinfeld generators are unchanged.} 
\be
e_1t_1=x_0^+,~~t_1^{-1}f_1=x_0^-,
~~e_0C=x_1^-,
~~C^{-1}f_0=x_{-1}^+,
~~t_1=t_1, 
~~t_0=Ct_1^{-1}.
\en

The subalgebra $U_A$ contains $(x^\pm_n)^{(r)}$ and $b_n$.
We shall work with the generating series
\be
x^{\pm}(z)=\sum_{n\in\Z}x^{\pm}_nz^{-n-1}.
\en

Set
\be
&&\Vt(\Lambda_i)
=K[b_{n}\mid n\in\Z_{<0}]\otimes
\bigl(\oplus_{m\in\Z} K e^{\Lambda_i+m\alpha_1}\bigr),
\\
&&
\Vt_A(\Lambda_i)=
A[b_{n}\mid n\in\Z_{<0}]\otimes
\bigl(\oplus_{m\in\Z} A e^{\Lambda_i+m\alpha_1}\bigr).
\en
Let $C$, $b_n$ ($n\in\Z\setminus\{0\}$) and $\partial$ act on an element 
$P\otimes e^{\beta}\in \Vt(\Lambda_i)$ as 
\be
&&C(P\otimes e^{\beta})=q(P\otimes e^{\beta}),
\\
&&b_n (P\otimes e^{\beta})=
\begin{cases}
(b_nP)\otimes e^{\beta}& (n<0),\\
\left[b_n,P\right]\otimes e^{\beta}& (n>0),\\
\end{cases}
\\
&&
\partial (P\otimes e^{\beta})=\langle h_1,\beta\rangle P\otimes e^{\beta}.
\en
We introduce the grading on $\Vt(\Lambda_i)$ by setting 
\be
\deg b_n=n,\quad 
\deg e^{\Lambda_i+m\alpha_1}=-m^2-im.
\en
Define the action of $U_q$ on $\Vt(\Lambda_i)$ by 
\bea
&&
x^+(z)\cdot u=
\exp\Bigl(-\sum_{n<0}\frac{b_n}{n}z^{-n}\Bigr)
\exp\Bigl(-\sum_{n>0}\frac{b_n}{n}(qz)^{-n}\Bigr)\,e^{\alpha_1}z^{\partial}\,u,
\label{vo1}
\\
&&x^-(z)\cdot u=
\exp\Bigl(\sum_{n<0}\frac{b_n}{n}(qz)^{-n}\Bigr)
\exp\Bigl(\sum_{n>0}\frac{b_n}{n}z^{-n}\Bigr)\,e^{-\alpha_1}z^{-\partial}\,u,
\label{vo2}
\\
&&
t_1\cdot u=q^{\partial }u,
~~
D\cdot u=q^{\deg u}\,u, 
\label{vo+}
\ena
where $u\in\Vt(\Lambda_i)$ is assumed to be homogeneous. 

\begin{prop}\cite{FJ,CJ}
For $i=0,1$, $\Vt(\Lambda_i)$ is a $U_q$-module 
isomorphic to $V(\Lambda_i)$, and $\Vt_A(\Lambda_i)$ 
is isomorphic to its $A$-form $V_A(\Lambda_i)$.  
\end{prop}

We will need also the realization of level $-1$ modules. 
Set 
\be
&&\Vt(-\Lambda_i)=K[b_{n}\mid n\in\Z_{>0}]\otimes
\bigl(\oplus_{m\in\Z} K e^{-\Lambda_i-m\alpha_1}\bigr),
\\
&&\Vt_A(-\Lambda_i)=A[b_{n}\mid n\in\Z_{>0}]\otimes
\bigl(\oplus_{m\in\Z} A e^{-\Lambda_i-m\alpha_1}\bigr),
\\
&&C(P\otimes e^{\beta})=q^{-1}(P\otimes e^{\beta}),
\\
&&b_n (P\otimes e^{\beta})=
\begin{cases}
(b_nP)\otimes e^{\beta}& (n>0),\\
\left[b_n,P\right]\otimes e^{\beta}& (n<0),\\
\end{cases},
\\
&&
\partial (P\otimes e^{\beta})=\langle h_1,\beta\rangle P\otimes e^{\beta},
\en
and 
\be
\deg b_n=n,\quad 
\deg e^{-\Lambda_i-m\alpha_1}=m^2+im.
\en
Let $U_q$ act on $u\in\Vt(-\Lambda_i)$ as 
\bea
&&x^+(z).u
=
\exp\Bigl(\sum_{n>0}\frac{b_n}{n}(q^{-1}z)^{-n}\Bigr)
\exp\Bigl(\sum_{n<0}\frac{b_n}{n}z^{-n}\Bigr)\,
z^{-\partial}e^{\alpha_1}\,u,
\label{vo3}
\\
&&x^-(z).u=
\exp\Bigl(-\sum_{n>0}\frac{b_n}{n}z^{-n}\Bigr)
\exp\Bigl(-\sum_{n<0}\frac{b_n}{n}(q^{-1}z)^{-n}\Bigr)
\,z^{\partial}e^{-\alpha_1}\,u,
\label{vo4}
\\
&&
t_1.u=q^{\partial }u,
~~
D.u=q^{\deg u}\,u. 
\label{vo-}
\ena

\begin{prop}
For $j=0,1$, $\Vt(-\Lambda_j)$ is a $U_q$-module 
isomorphic to $V(-\Lambda_j)$, and $\Vt_A(-\Lambda_j)$ 
is isomorphic to its $A$-form $V_A(-\Lambda_j)$.  
\end{prop}

For $i=0,1$ and $m\in\Z$ 
with $m\equiv i\bmod 2$, we set 
\be
v'_m=e^{\Lambda_i+(m-i)\alpha_1/2},
\quad
\vb'_m=e^{-\Lambda_i+(m+i)\alpha_1/2}.
\en
The following lemma shows that,  
up to a non-zero scalar multiple, they are extremal vectors
in $V(\Lambda_i)$ and $V(-\Lambda_i)$, respectively.  

\begin{lem}\label{lem:braid-action}
The following relations hold for $m\ge 0$. 
\be
&&(x_0^-)^{(m)}v'_m=(-q)^{m(m-1)/2}v'_{-m},
\\
&&(x_0^-)^{(m)}\vb'_m=(-q)^{-m(m-1)/2}\vb'_{-m},
\\
&&(x_{-1}^+)^{(m+1)}v'_{-m}=(-q)^{-m(m+1)/2}v'_{m+2},
\\
&&(x_{-1}^+)^{(m+1)}\vb'_{-m-2}=(-q)^{m(m+1)/2}\vb'_{m}. 
\en
\end{lem}
\begin{proof}
{}From \eqref{vo2}, it is straightforward to verify that
\be
&&x^-(z_1)\dots x^-(z_m)v'_m
\\
&&=
\prod_{j=1}^m z_j^{-m}\prod_{1\le j<k\le m}(z_j-z_k)(z_j-q^2z_k)
\times 
\exp\Bigl(-\sum_{n>0}\frac{b_{-n}}{n}\sum_{j=1}^m(qz_j)^{n}\Bigr)\,v'_{-m}.
\en
Hence $(x^-_0)^m v'_m$ is 
the coefficient of $(z_1\dots z_m)^0$ in 
\be
\prod_{j=1}^m z_j^{-m+1}\prod_{1\le j<k\le m}(z_j-z_k)(z_j-q^2z_k)
 \times v'_{-m}.
\en
We may replace this expression by 
its symmetrization with respect to $z_1,\dots,z_m$. 
Using 
\be
\sum_{\sigma\in S_m}\prod_{j<k}
\frac{z_{\sigma(j)}-q^2 z_{\sigma(k)}}
{z_{\sigma(j)}-z_{\sigma(k)}}=q^{m(m-1)/2}[m]!, 
\en
we find the first formula of Lemma. 
The other ones are obtained similarly. 
\end{proof}

\subsection{Proof of Proposition \ref{prop:support}}
\label{subsec:level1-2}

Proposition \ref{prop:support} will follow if we 
show the following formula: 
\begin{prop}\label{prop:xvm}
\bea
&&\Bigl(\sum_{n\ge 0}x_n^+z^{-n-1}\Bigr)(v'_i\otimes \vb'_m)
\label{eq:xvm1}
\\
&&\quad=q^{2m+i+2}z^{-m-2}\,
\exp\bigl(\sum_{k>0}\frac{b_k}{k}(q^{-1}z)^{-k}\bigr)\,
(v'_i\otimes \vb'_{m+2}),
\nn
\\
&&\Bigl(\sum_{n<0}x_n^+z^{-n-1}\Bigr)(v'_i\otimes \vb'_m)
\label{eq:xvm2}
\\
&&\quad =(-q^2)^{m+i+1}z^i
\,
\exp\bigl(\sum_{k>0}\frac{b_{-k}}{k}z^{k}\bigr)
(x^+_{-1})^{(m+2+i)}(x_0^-)^{(m+2+i)}\,(v'_i\otimes \vb'_{m+2}). 
\nn
\ena
\end{prop}
\begin{proof}
For the calculation, we use 
the following formulas proved in \cite{CP}. 
Let $U_q^{+,0}$ (resp., $U_q^{-,0}$) be the subalgebra generated by 
$e_i$ and $q^h$ 
(resp., $f_i$ and $q^h$) with $i\in I$, $h\in P^*$. 
Denote by $N^+_{\ge 0}$ (resp., $N^-_{> 0}$, 
$N^+_{<0}$, $N^-_{\le 0}$) the linear span of 
the elements $x^+_n$ ($n\ge 0$)
(resp., $x^-_n$ ($n>0$), $x^+_n$ ($n<0$), $x^-_n$ ($n\le 0$)).
Then 
\be
&&\Delta x^+_n\equiv
x_n^+\otimes C^n+\sum_{j=0}^nC^{2j}\varphi^+_{n-j}\otimes C^{n-j}x_j^+\\
&&\hspace{120pt}\bmod~~ U_q^{+,0}N^-_{>0}\otimes U_q^{+,0}(N^+_{\ge 0})^2
\quad(n\ge 0),
\\
&&\Delta x^+_{-n}\equiv
x_{-n}^+\otimes C^{-2n+1}
+\sum_{j=0}^{n-1}C^{2j-n}\varphi^-_{-j}\otimes C^{-2j}x_{-n+j}^+\\
&&\hspace{120pt}\bmod~~ U_q^{-,0}N^-_{\le0}\otimes U_q^{-,0}(N^+_{<0})^2
\quad(n>0),
\\
&&\Delta x^-_n\equiv
C^{2n-1}\otimes x_n^-+\sum_{j=0}^{n-1}C^{2j}x_{n-j}^-
\otimes C^{n-2j}\varphi_j^+\\
&&\hspace{120pt} \bmod~~ U_q^{+,0}(N^-_{>0})^2\otimes U_q^{+,0}N^+_{\ge 0}
\quad(n>0),
\\
&&\Delta x^-_{-n}\equiv
C^{-n}\otimes x_{-n}^-
+\sum_{j=0}^{n}C^{-n+j}x_{-j}^-\otimes C^{-2j}\varphi^-_{-n+j}\\
&&\hspace{120pt}\bmod~~ U_q^{-,0}(N^-_{\le0})^2\otimes U_q^{-,0}N^+_{<0}
\quad(n\ge 0),
\\
&&\Delta b_n\equiv
b_n\otimes C^n+C^{2n}\otimes b_n\\
&&\hspace{120pt}\bmod~~ U_q^{+,0}N^-_{>0}\otimes U_q^{+,0}N^+_{\ge 0}
\quad(n>0),
\\
&&\Delta b_{-n}\equiv
b_{-n}\otimes C^{-2n}+C^{-n}\otimes b_{-n}\\
&&\hspace{120pt}\bmod~~ U_q^{-,0}N^-_{\le 0}\otimes U_q^{-,0}N^+_{<0}
\quad (n>0). 
\en
{}From these we find that
\be
&&x^+_n(v'_i\otimes u)=
v'_i\otimes q^{2n+i}x^+_nu \quad (n\ge 0),
\\
&&b_n(v'_i\otimes u)=
v'_i\otimes q^{2n}b_nu \quad (n>0),  
\en
for $i=0,1$ and  $u\in V(-\Lambda_j)$. 
On the other hand, we have for $m\ge 0$ 
\be
\Bigl(\sum_{n\ge 0}x^+_nz^{-n-1}\Bigr)\vb'_m
&=&\Bigl(\sum_{n\in \Z}x^+_nz^{-n-1}\Bigr)\vb'_m
\\
&=&z^{-m-2}
\exp\Bigl(\sum_{k>0}\frac{b_k}{k}(q^{-1}z)^{-k}\Bigr)\vb'_{m+2}.  
\en
Eq. \eqref{eq:xvm1} follows from these relations. 
Calculating similarly, we find
\be
&&
x^+_{-n}(v'_i\otimes \vb'_m)=
(q^{2n-1}x^+_{-n}v'_i)\otimes \vb'_m \quad (n>0),
\\
&&b_{-n}(v'_i\otimes u)=
(q^{2n}b_{-n}v'_i)\otimes \vb'_m \quad (n>0),  
\en
and 
\be
\Bigl(\sum_{n< 0}x^+_nz^{-n-1}\Bigr)(v'_i\otimes\vb'_m)
=q^{2i+1}z^i
\exp\Bigl(\sum_{k>0}\frac{b_{-k}}{k}z^{k}\Bigr)
(v'_{i+2}\otimes\vb'_m).
\en
Acting with 
\be
&&\Delta((x^-_0)^{(r)})=\sum_{s=0}^rq^{s(r-s)}
(x^-_0)^{(s)}\otimes(x^-_0)^{(r-s)}t_1^{-s}, 
\\
&&\Delta((x_{-1}^+)^{(r)})=\sum_{s=0}^rq^{s(r-s)}
t_1^{-(r-s)}(x_{-1}^+)^{(s)}\otimes (x_{-1}^+)^{(r-s)}C^{-s},
\en
we obtain 
\be
&&
(x_0^-)^{(m+2+i)}(v'_i\otimes \vb'_{m+2})
=v'_{-i}\otimes (x_0)^{(m+2)}\vb'_{m+2},
\\
&&
(x_{-1}^+)^{(m+2+i)}(v'_{-i}\otimes (x_0)^{(m+2)}\vb'_{m+2})
=
(-1)^{m+i+1}q^{-2m-1}(v'_{i+2}\otimes \vb'_m),
\en
where we have used 
\be
&&(x_{-1}^+)^{(m+1)}(x_0^-)^{(m+2)}\vb'_{m+2}
=(-q)^{-m-1}\vb'_m,
\\
&&(x_{-1}^+)^{(i+1)}v'_{-i}=(-q)^{-i}v'_{i+2},
\en
which follow from Lemma \ref{lem:braid-action}.
Combining these relations we arrive at
\eqref{eq:xvm2}.
\end{proof}

The proof of Proposition \ref{prop:support} 
is now complete.
\subsection{Intertwining operators}\label{subsec:intertwiners}

In this section we show Proposition \ref{prop:FtoV}. 

For integrable modules $M,N$, define the completion of $M\otimes N$
by 
$(M\otimes N)^{\wedge}=
\sum_{\mu,\nu}\prod_{\xi\in Q_+}M_{\mu+\xi}\otimes N_{\nu-\xi}$, 
where $Q_+=\Z_{\ge 0}\alpha_0+\Z_{\ge 0}\alpha_1 $. 
Similarly define $(M_1\otimes\dots\otimes M_p)^\wedge
=((M_1\otimes\dots\otimes M_{p-1})^\wedge\otimes M_p)^\wedge$.
{}From the definition we have 
\be 
(V_{z}^{\otimes n})^{\wedge}= 
V^{\otimes n} \otimes_K
K[[z_1/z_2,\ldots,z_{n-1}/z_n]][z_1^{\pm 1},\dots,z_n^{\pm 1}]. 
\en 
Consider the intertwiner of $U_q$-modules of the form 
\bea
&&\Psit(z):
V(\Lambda_i)\longrightarrow 
\bigl(V_z\otimes V(\Lambda_{1-i})\bigr)^\wedge.  
\label{typeII}
\ena
We write 
$\Psit(z)u=\sum_{\epsilon,n} 
z^{-n}v_\epsilon\otimes\Psit_{\epsilon,n} u$ 
for $u\in V(\Lambda_i)$, and set
$\Psit_\epsilon(z)=\sum_{n\in\Z}\Psit_{\epsilon,n}z^{-n}$. 
We choose  
the normalization $\Psit_{-,0}v_0=v_1$, $\Psit_{+,0}v_1=v_0$. 
For more details concerning $\Psit(z)$,  
we refer to \cite{JM}, Chapter 6. 

Denote by
\bea\label{PAIRING}
\langle~,~\rangle~:~ V(\Lambda_i)\otimes V(-\Lambda_i)\rightarrow K
\ena
a $U_q$-linear mapping normalized as $\langle v_i,\vb_{-i}\rangle=1$. 

Iterating \eqref{typeII}, we obtain a $U'_q$-linear map 
\bea
&&\psi_n:
V(\Lambda_i)\otimes V(-\Lambda_j)
\longrightarrow (V_z^{\otimes n})^\wedge
\qquad (n\equiv i-j\bmod 2),
\nn \\
&&\psi_n(u\otimes v)
=\rho^{(i,j)}_n(z_1,\dots,z_n)
\sum_{\epsilon_1,\dots,\epsilon_n}
\langle\Psit_{\epsilon_n}(z_n)\dots
\Psit_{\epsilon_1}(z_1)u,v\rangle
v_{\epsilon_1}\otimes\dots\otimes v_{\epsilon_n},
\nn 
\ena 
where 
\bea 
\rho^{(i,j)}_n(z_1,\dots,z_n)
&=&
(-q)^{-l(l-j)}
\prod_{\substack{1\le k\le n\\\mbox{\tiny $k$:odd}}}z_k^{-(k-1)/2}
\prod_{\substack{1\le k\le n\\\mbox{\tiny $k$:even}}} z_k^{-k/2+i} 
\label{def:rho} \\ 
&\times& 
\prod_{1\le k<k'\le n}
\frac{(q^2z_k/z_{k'};q^4)_\infty}{(z_k/z_{k'};q^4)_\infty}. \nn
\ena
Here we have set $n=2l+i-j$, 
and $(z;p)_\infty=\prod_{n\ge0}(1-p^nz)$
(note that $(z;q^4)_\infty^{\pm 1}\in A[[z]]$). 
We have included the scalar factor 
$\rho^{(i,j)}_n(z_1,\dots,z_n)$, so that 
the normalization condition
\bea
&&\psi_n(v_i\otimes \vb_{n-i})=
\overbrace{v_+\otimes \dots\otimes v_+}^{{\mbox{\tiny $n$ times}}}
\label{normal}
\ena
holds. 
This has the effect of shifting degrees as
\bea
\deg (u\otimes v)
=\deg \psi_n(u\otimes v)+\frac{(n-i)^2-j}{4}.
\label{shift}
\ena
Proposition \ref{prop:FtoV} is a consequence of \eqref{FntoVn} below.
(Note that $\tpsi_{n}$ in Proposition \ref{prop:FtoV} is the restriction of $\psi_{n}$ to $F_{n}^{(i, j)}$.) 
\begin{prop}\label{prop:intertwiner}
We have
\bea
&&\psi_n(F^{(i,j)}_n)=V(n(\Lambda_1-\Lambda_0)),\label{FntoVn}\\
&&\psi_n(F^{(i,j)}_{n+2})=0,\label{n+1}\\
&&\psi_n\bigl(V_A(\Lambda_i)\otimes V_A(-\Lambda_j)\bigr)
\subset (V_{A,z}^{\otimes n})^\wedge.\label{VAtoVAaf}
\ena
\end{prop}
\begin{proof}
Since $V\bigl(n(\Lambda_0-\Lambda_1)\bigr)$
is generated by $v_+^{\otimes n}$, 
the assertion \eqref{FntoVn} is clear from 
\eqref{normal}. 
The assertion (\ref{n+1}) follows from the fact that
$\psi_n(v_0\otimes \vb_{n+2})=0$,
which is obvious from consideration of weights.
To see \eqref{VAtoVAaf}, because of the cyclicity of 
$V_A(\Lambda_i)\otimes V_A(-\Lambda_j)$, it suffices to show that the 
vector $\psi_n(v_i\otimes \vb_{-j})$ belongs to 
$(V_{A,z}^{\otimes n})^\wedge$. 
Set 
\bea
&&\rho^{(i,j)}_n(z_1,\dots,z_n)
\langle\Psit_{\epsilon_n}(z_n)\cdots
\Psit_{\epsilon_1}(z_1)v_i,\vb_{-j}\rangle
=
\frac{a_{\epsilon_1,\dots,\epsilon_n}(z_1,\dots,z_n)}
{\prod_{1\le r<s\le n}(1-q^{-2}z_r/z_s)}.
\label{mat-elt}
\ena
As before, let $n=2l+i-j$. 
{}From \cite{JM}, eq.(9.8) and p.116,  we have
\bea
a_{\underbrace{-\dots-}_{l}\underbrace{+\dots+}_{n-l}}(z_1,\dots,z_n)
&=&(-q)^{l(l-1)/2}
\prod_{1\le r\le l}z_r^i
\prod_{l+1\le r\le n}z_r^{-l}
\nn\\
&\times&
\prod_{1\le r<s\le l}(1-q^{-2}z_r/z_s)
\prod_{l+1\le r<s\le n}(1-q^{-2}z_r/z_s),
\label{mat-elt2}\\
a_{\dots,\pm,\mp,\dots}(\dots,z_k,z_{k+1},\dots)
&=&
\frac{(z_k/z_{k+1}-q^2)(z_k/z_{k+1})}{q(1-z_k/z_{k+1})}
a_{\dots,\mp,\pm,\dots}(\dots,z_{k+1},z_{k},\dots)
\nn\\
&&
-\frac{(1-q^2)(z_k/z_{k+1})^t}{q(1-z_k/z_{k+1})}
a_{\dots,\mp,\pm,\dots}(\dots,z_{k},z_{k+1},\dots). 
\label{mat-elt3}
\ena
Here $t=0$ for the upper sign and $1$ for the lower sign. 
In the right hand side, there is no pole at $1-z_k/z_{k+1}=0$. 
It follows that 
$a_{\epsilon_1,\dots,\epsilon_n}(z_1,\dots,z_n)$ are 
Laurent polynomials in $q,z_1,\dots,z_n$.
Expanding the right hand side of \eqref{mat-elt}
into a Laurent series in $z_1,\dots,z_n$, 
we see that all coefficients are Laurent polynomials in $q$. 
\end{proof}
\section{Functional model}\label{sec:functional}
\subsection{The space $\F_{n}$}
\label{subsec:TV}
In this and the next subsections, 
we introduce various spaces of 
polynomials used to give a realization of the tensor product 
$V(\Lambda_i)\otimes V(-\Lambda_j)$. 
As before, we set $K=\C(q),A=\C[q,q^{-1}]$. 

For $0\le l\le n$, let $\F_{n,l}$ be the space
of polynomials $P(X_1,\ldots,X_l)$ with coefficients in 
$K[z_1^{\pm 1},\ldots,z_n^{\pm 1}]$, satisfying the following
conditions:
\bea
&&\mbox{$P$ is symmetric in $X_1,\ldots,X_l$},
\label{PsymX}\\
&&\deg_{X_i}P\le n-1,
\label{degn}\\
&&P\vert_{X_1=q^{-2}X_2=z_k^{-1}}=0\quad\mbox{for $1\le k\le n$ and $l\ge 2$.}
\label{XXZ}
\ena
We set 
\be
&&\F_{A,n,l}:=\F_{n,l}\cap 
A[z_1^{\pm 1},\ldots,z_n^{\pm 1}][X_1,\ldots,X_l],
\\
&&
\F_n:=\oplus_{l=0}^n\F_{n,l},
\quad
\F_{A,n}:=\oplus_{l=0}^n\F_{A,n,l}.
\en

For each subset $M\subset \{1,\ldots,n\}$ with $\sharp M=l$, 
define $w^{(n)}_M(X_1,\ldots,X_l)\in \F_{A,n,l}$ by 
\bea
&&w^{(n)}_M(X_1,\ldots,X_l)
:=\Sym \Bigl(G^{(n)}_{m_1}(X_1)\cdots G^{(n)}_{m_l}(X_l)
\prod_{1\le k<k'\le l}\frac{q^{-1}X_k-qX_{k'}}{X_k-X_{k'}}\Bigr).
\label{wM}
\ena
Here $M=\{m_1,\ldots,m_l\}$ ($1\le m_1<\cdots<m_l\le n$), 
${\rm Sym}$ stands for the symmetrization 
\be
(\Sym f)(X_1,\ldots,X_l):=\sum_{\sigma\in \mathfrak S_l} f(X_{\sigma(1)},
\ldots,X_{\sigma(l)}),
\en
and 
\be
&&G^{(n)}_m(X):=q^{m-n}\prod_{k=1}^{m-1}(1-q^{-2}z_kX)\prod_{k=m+1}^n(1-z_kX). 
\en
We will write \eqref{wM} also as  
$w^{(n)}_{\epsilon_1,\ldots,\epsilon_n}(X_1,\ldots,X_l)
=w^{(n)}_{\epsilon_1,\ldots,\epsilon_n}(X_1,\ldots,X_l|z_1,\ldots,z_n)$, 
where $\epsilon_1,\ldots,\epsilon_n\in\{+,-\}$ 
are related to $M$ via $M=\{j\mid \epsilon_j=-\}$. 
The polynomials \eqref{wM} arise naturally 
in the framework of the quantum inverse scattering method
(see (C.1) and (C.2) in \cite{JMT}).  
We have the transformation property 
\bea
&&
w^{(n)}_{\ldots,\epsilon'_{j+1},\epsilon'_j,\ldots}
(X_1,\ldots,X_l|\ldots,z_{j+1},z_j,\ldots)
\nn
\\
&&\quad=\sum_{\varepsilon_j,\varepsilon_{j+1}}
w^{(n)}_{\ldots,\epsilon_{j},\epsilon_{j+1},\ldots}
(X_1,\ldots,X_l|\ldots,z_{j},z_{j+1},\ldots)
\bigl(R(z_j/z_{j+1})^{-1}\bigr)
_{\epsilon_j,\epsilon_{j+1};\epsilon'_j,\epsilon'_{j+1}},
\label{exchange}
\ena
where 
\be
R(z)=\begin{pmatrix}
1 & & & \\
 & \displaystyle{\frac{(1-z)q}{1-q^2z}}&
\displaystyle{\frac{1-q^2}{1-q^2z}}& \\
 & 
\displaystyle{\frac{(1-q^2)z}{1-q^2z}}&
\displaystyle{\frac{(1-z)}{1-q^2z}}& \\
 & & & 1\\
\end{pmatrix}.
\en
We assign a degree to an element $P\in\F_n$ by setting 
\bea
\deg X_p=-1,\quad \deg z_j=1.
\label{degree-assign}
\ena

Define a completion of $\widehat{\F}_n=\oplus_{l=0}^n\widehat{\F}_{n,l}$ by
\be
\widehat{\F}_{n,l}:=
\F_{n,l}\bigotimes_{K[z_1^{\pm 1},\ldots,z_n^{\pm 1}]}
K[[z_1/z_2,\ldots,z_{n-1}/z_n]][z_1^{\pm 1},\ldots,z_n^{\pm 1}].
\en
Similarly, we define the completion $\widehat{\F}_{A,n}$.

We have 
\begin{lem}\label{lem:basis}\cite{NPT}
For each $l$ $(0\leq l\leq n)$, the set of polynomials
$w^{(n)}_{M}(X_1,\dots,X_l)$ with $\# M=l$ 
constitutes a free basis of $\widehat\F_{A,n,l}$ over
\be 
\widehat A_n=A[[z_1/z_2,\ldots,z_{n-1}/z_n]][z_1^{\pm 1},\ldots,z_n^{\pm 1}]. 
\en 
\end{lem}
\begin{proof}
Let $\Pc_{n,l}$ be the set of partitions
$\lambda=(\lambda_1,\ldots,\lambda_l)$ satisfying $\lambda_1\leq n-1$.
Consider an element of $\widehat A_n[X_1,\ldots,X_l]$ of the form
\bea\label{MONO}
P_{n,l}=\sum_{\lambda\in\Pc_{n,l}}c_\lambda m_\lambda,
\ena
where $m_\lambda$ is the monomial symmetric polynomial corresponding to
$\lambda$. The condition (\ref{XXZ})
is equivalent to a set of linear relations
for the coefficients $c_\lambda\in \widehat A_n$.
The linear relations are defined over $A[z_1^{\pm 1},\ldots,z_n^{\pm 1}]$.

Let  $d_{n,l}$ be the dimension of the space of solutions
$(c_\lambda)_{\lambda\in P_{n,l}}$ where $c_\lambda$ belongs to the quotient
field of $\widehat A_n$.
If we specialize the relation (\ref{XXZ}) to $q=\sqrt{-1}$,
it is equivalent to
\bea\label{eq:Pnl}
&&P_{n,l}|_{X_1=-X_2}=0.
\ena
This is because $P_{n,l}\bigl|_{X_1=X,X_2=-X}$ is a polynomial in $X$
of degree at most $2(n-1)$ and has $2n$ zeroes $\pm z_k^{-1}$. 
Therefore, we have $d_{n,l}\leq\binom nl$.

Note that $G^{(n)}_m(z_j^{-1})=0$ if $m<j$, and
$G^{(n)}_m(z_m^{-1})$ is an invertible element in
$\widehat{A}_{n}$. 

For subsets $M=\{m_1<\cdots<m_l\}$,
$J=\{j_1<\cdots<j_l\}$, we write 
$M\le J$ if and only if $m_a\le j_a$ for all $a$. 
Then $\le$ is a partial ordering. 
By induction on $l$ one can show that
\be
w^{(n)}_J(z_{m_1}^{-1},\ldots,z_{m_n}^{-1})=0
\en
unless $M\leq J$, and
$w^{(n)}_J(z_{j_1}^{-1},\ldots,z_{j_n}^{-1})$ is invertible in $\widehat A_n$.
Using this triangularity one can show that if
$\sum_{M,\sharp(M)=l}c_Mw^{(n)}_M=0$, then $c_M=0$ for all $M$. Therefore,
we have $d_{n,l}\geq\binom nl$, and $d_{n,l}=\binom nl$.
In conclusion, we proved that $w^{(n)}_M$ such that $\sharp(M)=l$ constitute
a basis 
of the vector space of solutions to (\ref{XXZ}) over the quotient field of $\widehat A_n$.

Let $f\in\F_{A,n,l}$ be written as
\be
f=\sum_Mc_Mw^{(n)}_M.
\en
We show that $c_M\in\widehat A_n$.
Suppose that $J\subset\{1,\ldots,n\}$ is a maximum element in
$\{M;c_M\not=0\}$. We have $f(z_{j_1}^{-1},\ldots,z_{j_l}^{-1})=
c_Jw^{(n)}_J(z_{j_1}^{-1},\ldots,z_{j_n}^{-1})$,
$f(z_{j_1}^{-1},\ldots,z_{j_l}^{-1})\in\widehat A_n$, and moreover,
$w^{(n)}_J(z_{j_1}^{-1},\ldots,z_{j_n}^{-1})$ is invertible in $\widehat A_n$.
Therefore, we have $c_J\in\widehat A_n$. Applying the same argument to
$f-c_Jw^{(n)}_J\in\widehat A_n[X_1,\ldots,X_l]$, and proceeding inductively,
we see that $c_M\in\widehat A_n$
for all $M$.
\end{proof}

The following result can be extracted from \cite{TV}
and Propositions C.1, C.2 in \cite{JMT}. 
There the symmetric group $\mathfrak S_n$ acts on $\F_n$ by the permutation of variables
$z_1,\ldots,z_n$.

\begin{prop}\label{prop:TV}
$(\hbox{\rm i})$
The space $\F_n$ has a structure of a
$U_q$-module of level $0$. 
The action of $U_q'$ commutes with the multiplication 
by elements of $K[z_1^{\pm 1},\ldots,z_n^{\pm 1}]$, 
and $D P=q^d P$, where $d=\deg P$ is given by \eqref{degree-assign},
and $t_1$ acts on $\F_{n,l}$ by the multiplication of $q^{n-2l}$.
The action preserves the $\mathfrak S_n$-invariant subspace$:$
\bea\label{INV}
U_q\F_n^{\mathfrak S_n}\subset\F_n^{\mathfrak S_n}.
\ena

\noindent$(\hbox{\rm ii})$
There exists an injective morphism of $U_q$-modules
given by
\bea
&&\Cc_n:V_z^{\otimes n} \longrightarrow \F_n,
\qquad
v_{\epsilon_1}\otimes\cdots\otimes v_{\epsilon_n}
\mapsto 
w^{(n)}_{\epsilon_1,\ldots,\epsilon_n}(X_1,\ldots,X_l).
\label{TVmap}
\ena
Moreover $\Cc_n$ is $K[z_1^{\pm 1},\ldots,z_n^{\pm 1}]$-linear.

\noindent$(\hbox{\rm iii})$
We have 
\be
U_A\F_{A,n}\subset \F_{A,n},
\quad \Cc_n(V_{A,z}^{\otimes n})\subset \F_{A,n}.
\en

\end{prop}
Explicit formulas for the action on $\F_n$ can be found 
in \cite{JMT}, Proposition C.1. 

By Lemma \ref{lem:basis}, the morphism $\Cc_n$ can be extended to an
isomorphism between the completions.
\be
\widehat{\Cc}_n :(V_z^{\otimes n})^{\widehat{}}
\isoto\widehat{\F}_n,
\en

Denote by $\varphi_n:=\widehat{\Cc}_n\circ\psi_n$ the composition
\be
\varphi_n\colon V(\Lambda_i)\otimes V(-\Lambda_j)
\To[\psi_n]
(V_z^{\otimes n})^{\widehat{}}
\To[\Cc_n]\widehat{\F}_n
\qquad (n\equiv i-j\bmod 2).
\en
Then $\varphi_n$ is $U'_q$-linear. 
For a weight vector 
$u\otimes v \in V(\Lambda_i)\otimes V(-\Lambda_j)$, 
let $m=\langle h_1,\wt(u\otimes v)\rangle$.   
Then 
$\varphi_{n}(u\otimes v)\in\oplus_{l=0}^n\widehat\F_{n,l}$ 
has the only non-zero component for $l=(n-m)/2$. 
We denote it by 
$P^{u\otimes v}_{n,l}=
P^{u\otimes v}_{n,l}(X_1,\ldots,X_l|z_1,\ldots,z_n)$. 
Explicitly we have 
\bea
P^{u\otimes v}_{n,l}
&=&\rho^{(i,j)}_n(z_1,\ldots,z_n) \label{Puv} \\ 
&\times& 
\sum_{\epsilon_1,\ldots,\epsilon_n}
\langle\Psit_{\epsilon_n}(z_n)\cdots
\Psit_{\epsilon_1}(z_1)u,v\rangle
w^{(n)}_{\epsilon_1,\ldots,\epsilon_n}
(X_1,\ldots,X_l|z_1,\ldots,z_n). \nn 
\ena
Note also that 
\be
\deg (u\otimes v) =\frac{(n-i)^2-j}{4}+
\deg P^{u\otimes v}_{n,l}.
\en

Let us consider the special case $v_i\otimes\vb_{-j}$. 
\begin{prop}\label{i}
For $i,j\in\{0,1\}$ and $n\ge 0$ with $n=2l+i-j$, we have
\bea
\varphi_n(v_i\otimes\vb_{-j})
=(-)^{l(l-1)/2}q^{\gamma}
\prod_{a=1}^lX_a^{1-j}\prod_{1\le a\neq b\le l}(X_a-q^{-2}X_b),  
\label{top}
\ena
where $\gamma=l(l-1)/2-2l(1-j)-l(n-1)$.
\end{prop}
\begin{proof}
The right hand side of \eqref{top} belongs to $\F_n$.
{}From Lemma \ref{lem:basis}, we can therefore write
\bea
(-)^{l(l-1)/2}q^{\gamma}
\prod_{a=1}^lX_a^{1-j}\prod_{1\le a\neq b\le l}(X_a-q^{-2}X_b) 
=\sum_{\epsilon_1,\ldots,\epsilon_n}
c_{\epsilon_1,\ldots,\epsilon_n}
w^{(n)}_{\epsilon_1,\ldots,\epsilon_n},
\label{tobedet}
\ena
where $c_{\epsilon_1,\ldots,\epsilon_n}$ 
are some rational functions in $q,z_1,\ldots,z_n$. 
We are to show the relation 
\bea
c_{\epsilon_1,\ldots,\epsilon_n}=
\rho^{(i,j)}_n(z_1,\ldots,z_n)
\langle \Psit_{\epsilon_n}(z_n)
\cdots\Psit_{\epsilon_1}(z_1)v_i,\vb_{-j}\rangle.
\label{matrix-element}
\ena
Specializing $X_1=q^2z_1^{-1},\ldots, X_l=q^2z_l^{-1}$
in \eqref{tobedet}, we find that only the term with 
$\epsilon_1=\cdots=\epsilon_l=-$, $\epsilon_{l+1}=\cdots=\epsilon_n=+$ 
contributes. 
Comparing the result with the matrix element \eqref{mat-elt2}, 
we obtain \eqref{matrix-element} for 
$\epsilon_1=\cdots=\epsilon_l=-$. 
The general case follows from this, since
both sides of \eqref{matrix-element} 
share the same transformation property 
under the exchange of $(z_k,\epsilon_k)$ and 
$(z_{k+1},\epsilon_{k+1})$.  
\end{proof}

In general, the image of $\varphi_n$ has the following properties. 
\begin{prop}\label{prop:generic-q}
$(\hbox{\rm i})$ 
\be
\varphi_n(V(\Lambda_i)\otimes V(-\Lambda_j))\subset \F_n^{\mathfrak S_n},
\quad
\varphi_n(V_A(\Lambda_i)\otimes V_A(-\Lambda_j))\subset \F_{A,n}^{\mathfrak S_n}.
\en

\noindent$(\hbox{\rm ii})$ We have 
\bea
\varphi_n(F^{(i,j)}_{n+2})=0.
\label{phiFn}
\ena
The induced map
$\overline\varphi_n:F^{(i,j)}_n/F^{(i,j)}_{n+2}\rightarrow\F_n^{\mathfrak S_n}$
is injective.

\noindent$(\hbox{\rm iii})$
Let $u\otimes v\in V(\Lambda_i)\otimes V(-\Lambda_j)$ and  
$P^{u\otimes v}_{n,l}$ be as above. We have 
\bea
&&P^{{u\otimes v}}_{n+2,l+1}
(X_1,\ldots,X_l,z^{-1}|z_1,\ldots,z_n,z,q^{2}z)
\label{LINK}\\
&&\quad=q^{\nu}z^{-n-1+i}
\prod_{a=1}^l(1-q^{-2}zX_a)(1-q^{2}zX_a)
\times P^{{u\otimes v}}_{n,l}(X_1,\ldots,X_l|z_1,\ldots,z_n).
\nonumber
\ena
Here $n=2l+i-j$, $\nu=-5l+i+j-2$ for $n$ even, 
and $\nu=-5l+i+j-1$ otherwise. 
\end{prop}

\begin{proof}
By Proposition \ref{i}, we have 
$\varphi_n(v_i\otimes\vb_{-j})\in \F_{A,n}^{\mathfrak S_n}$.
Since 
$V(\Lambda_i)\otimes V(-\Lambda_j)
=U_q(v_i\otimes \vb_{-j})$
and 
$V_A(\Lambda_i)\otimes V_A(-\Lambda_j)
=U_A(v_i\otimes \vb_{-j})$, assertion (i) follows from (\ref{INV}).

We have \eqref{phiFn} by \eqref{n+1}.
Moreover the composition
\be
V(n(\Lambda_1-\Lambda_0))\overset{\sim}{\longrightarrow}
F^{(i,j)}_n/F^{(i,j)}_{n+2}\overset{\overline\varphi_n}
{\longrightarrow}\F_n^{\mathfrak S_n}
\en
coincides with $\Cc_n$ which is injective. 
Hence $\overline\varphi_n$ is also injective, 
and we have (ii). 

Let us show (\ref{LINK}).
It is enough to verify it assuming that 
$q$ is a complex number with $|q|<1$. 
In this case, the intertwiners satisfy the properties
(\cite{JM}, eq.(6.40), with $v\in V(\Lambda_j)$):
\be
&&\Psit_\epsilon(z')\Psit_\epsilon(z)v=O(1)
\qquad (z'\to q^2z), \\
&&\frac{(q^2z/z';q^4)_{\infty}}
{(q^4z/z';q^4)_{\infty}}
\bigl(\Psit_+(z')\Psit_-(z)-q^{-1}\Psit_-(z')\Psit_+(z)\bigr)v\\
&&\hspace{50pt}=(-q)^{-j}v+O\bigl((z'-q^2z)\bigr)
\qquad (z'\to q^2z). 
\en
It is easy to verify the relations 
\be
&&w^{(n+2)}_{\epsilon_1,\ldots,\epsilon_{n},\epsilon,\epsilon}
(X_1,\ldots,X_l,z^{-1}|z_1,\ldots,z_n,z,q^2z)
=0,
\\
&&
-qw^{(n+2)}_{\epsilon_1,\ldots,\epsilon_{n},+\,-}
(X_1,\ldots,X_l,z^{-1}|z_1,\ldots,z_n,z,q^2z)
\\
&&\quad=
w^{(n+2)}_{\epsilon_1,\ldots,\epsilon_{n},-\,+}
(X_1,\ldots,X_l,z^{-1}|z_1,\ldots,z_n,z,q^2z)
\\
&&
\quad=q^{-l-1}(1-q^2)\prod_{k=1}^n(1-q^{-2}z_k/z)
\prod_{a=1}^l(1-q^{-2}zX_a)(1-q^{2}zX_a)
\\
&&\qquad\times 
w^{(n)}_{\epsilon_1,\ldots,\epsilon_{n}}
(X_1,\ldots,X_l|z_1,\ldots,z_n).
\en
Specializing \eqref{Puv} to $X_{l+1}=z^{-1}, z_{n+1}=z, z_{n+2}=q^{2}z$, 
and using these relations along with \eqref{def:rho}, 
we obtain (\ref{LINK}). 
\end{proof}

\subsection{The spaces $\widehat{\Zc}^{(i,j)}$ and 
$W_n$}

Proposition \ref{prop:generic-q} motivates us
to consider the subspace $\widehat{\Zc}^{(i,j)}[m]\subset
\prod\F_{n,l}^{\mathfrak S_n}$,  
consisting of all sequences of polynomials 
$(P_{n,l})_{\substack{n\ge 0\\ n-2l=m}}$ such that 
\bea
&&P_{n+2,l+1}\hbox{ and } P_{n,l}\hbox{ are related 
by (\ref{LINK})}.
\label{link2}
\ena

We set 
\be
&&\widehat{\Zc}_A^{(i,j)}[m]=
\widehat{\Zc}^{(i,j)}[m]\cap 
\prod_{\substack{n\ge 0\\ n-2l=m}}\F_{A,n,l}^{\mathfrak S_n}, 
\\
&&\widehat{\Zc}^{(i,j)}=
\bigoplus_{\substack{m\in\Z\\ m\equiv i-j\bmod 2}}
\widehat{\Zc}^{(i,j)}[m],
\\
&&\widehat{\Zc}_A^{(i,j)}
=\bigoplus_{\substack{m\in\Z\\m\equiv i-j\bmod 2}}
\widehat{\Zc}_A^{(i,j)}[m]. 
\en
The space $\widehat{\Zc}^{(i,j)}$ is a $U_q$-module
by the componentwise action. 

Define also a filtration of $\widehat{\Zc}^{(i,j)}$ by setting 
\bea
&&\widehat{\Zc}^{(i,j)}_n
=\bigoplus_{\substack{m\in\Z\\m\equiv i-j\bmod 2}}
\widehat{\Zc}^{(i,j)}_n[m],\\
&&
\begin{array}{l}
\widehat{\Zc}^{(i,j)}_n[m]
=\{(P_{n',l})_{\substack{n'\ge 0\\n'-2l=m}}
\in\widehat{\Zc}^{(i,j)}[m]\mid\\
\hspace{100pt}P_{n',l}=0\hbox{ for all $n',l$ such that $n'<n$}\}.
\end{array}
\ena
{}From Proposition \ref{prop:generic-q}, we have a morphism 
of $U'_q$-modules
\be
\varphi=\prod_{n}\varphi_n~:~
V(\Lambda_i)\otimes V(-\Lambda_j)
\longrightarrow \widehat{\Zc}^{(i,j)}, 
\en
which satisfies $\varphi\bigl(F^{(i,j)}_n\bigr)\subset 
\widehat{\Zc}_n^{(i,j)}$.  
By Proposition \ref{prop:generic-q} (ii) and 
Proposition \ref{prop:intersection}, $\varphi$ is injective. 
Our goal is to show that,
after an appropriate completion of
$V(\Lambda_i)\otimes V(-\Lambda_j)$, the mapping $\varphi$
becomes an isomorphism
(see Theorem \ref{thm:functional-realization} below). 

For that purpose, 
let us introduce the following spaces.  
Let $W_{n,l}$ denote the subspace of $\F_{n,l}^{\mathfrak S_n}$ consisting of 
elements $P$ such that
\be
&&P\vert_{X_1=z_1{}^{-1}=(q^{-2}z_2)^{-1}}=0\quad\mbox{for $n\ge2$ and $l\ge1$.}
\en
Set further
\be
&&W_{A,n,l}:=W_{n,l}\cap \F_{A,n,l},
\\
&&
W^{\ge 0}_{n,l}:=W_{n,l}\cap K[z_1,\ldots,z_n][X_1,\ldots,X_l],
\\
&&
W^{\ge 0}_{A,n,l}:=W_{n,l}\cap A[z_1,\ldots,z_n][X_1,\ldots,X_l].
\en
We set $W_n:=\oplus_{l=0}^nW_{n,l}$, and similarly for
$W_{A,n}$, $W^{\ge 0}_n$, $W_{A,n}^{\ge 0}$. 
{}From the explicit action (see \cite{JMT}, Proposition C.1)
we have $U_qW_n\subset W_n$, $U_AW_{A,n}\subset W_{A,n}$
and $U_q^{\ge 0}W^{\ge 0}_n\subset W^{\ge 0}_n$. 

Consider the isomorphisms of $U'_q$-modules
\bea\label{B}
&&V\bigl(n(\Lambda_1-\Lambda_0)\bigr)
\overset{\sim}{\longrightarrow} 
F^{(i,j)}_n/F^{(i,j)}_{n+2}
\overset{\sim}{\longrightarrow} 
\varphi(F^{(i,j)}_n)/\varphi(F^{(i,j)}_{n+2}). 
\ena
The first map given by Proposition \ref{prop:FtoV} 
shifts the degree by $+s$, 
and the second map $\varphi$ by $-s$, where $s=((n-i)^2-j)/4$. 
Hence the composition is $U_q$-linear. 
By Proposition \ref{prop:generic-q} (ii), there are also
injective canonical maps
\bea\label{A}
\varphi(F^{(i,j)}_n)/\varphi(F^{(i,j)}_{n+2})
\mono
\widehat{\Zc}^{(i,j)}_n/\widehat{\Zc}^{(i,j)}_{n+2}
\mono W_n. 
\ena
The composition of (\ref{B}) and (\ref{A}) 
\bea
V(n(\Lambda_1-\Lambda_0))\longrightarrow W_n
\label{C}
\ena
coincides with the restriction
to $V(n(\Lambda_1-\Lambda_0))
\simeq U_qv_+^{\otimes n}\subset V_z^{\otimes n}$
of the map $\Cc_n$ defined in (\ref{TVmap}).

\begin{prop}\label{prop:W-is-cyclic}
Let $1_n$ be the unit of $\F_{n,0}$. We have
\be
W_n=U_q1_n, 
\quad W^{\ge 0}_n=U^{\ge 0}_q1_n.
\en
\end{prop}
We defer the proof of Proposition \ref{prop:W-is-cyclic}
to the next subsection. 

\begin{thm}\label{thm:V=W}
The morphisms \eqref{A}, \eqref{C} are isomorphisms. 
\end{thm}
\begin{proof}
The map \eqref{C} is injective because so is $\Cc_n$.
Proposition \ref{prop:W-is-cyclic} shows that 
it is also surjective. 
\end{proof}

Define the completed tensor product with respect to the filtration 
$\{F_n^{(i,j)}\}$
\be
V(\Lambda_i)\ftens V(-\Lambda_j)
=\varprojlim
V(\Lambda_i)\otimes V(-\Lambda_j)/F_n^{(i,j)}.
\en
By Proposition \ref{prop:intersection}, 
the map $V(\Lambda_i)\otimes V(-\Lambda_j)\to V(\Lambda_i)\ftens V(-\Lambda_j)$
is injective.
Clearly 
$\varprojlim\widehat{\Zc}^{(i,j)}/\widehat{\Zc}^{(i,j)}_n
=\widehat{\Zc}^{(i,j)}$.  
Theorem \ref{thm:V=W} implies that 
\be
V(\Lambda_i)\ftens V(-\Lambda_j)
\longrightarrow \widehat{\Zc}^{(i,j)} 
\en
is an isomorphism.
Hence we arrive at the following result, 
which provides a `functional realization' of the 
(completed) tensor product of level $1$ and level $-1$
integrable modules.

\begin{thm}\label{thm:functional-realization}
We have an isomorphism 
\be
&&V(\Lambda_i)\ftens V(-\Lambda_j)\isoto\widehat{\Zc}^{(i,j)}.
\en
\end{thm}

\subsection{Relation with $\infty$-cycles}

In \cite{JMMT}, sequences similar to the elements of 
$\widehat{\Zc}^{(i,j)}$ have been considered under the name 
`$\infty$-cycles'. 
The latter are closely related to the specialization of 
the former at $q=\sqrt{-1}$. 
The aim of this subsection is to clarify the 
connection between these objects. 

First we recall the definitions given in \cite{JMT,JMMT}
\footnote{The spaces 
$\F^{\rm skew}_{\C,n,l}$, $W^{\rm skew}_{\C,n,l}$,
$W^{{\rm skew}\ge 0}_{\C,n,l}$ and
$\widehat{\Zc}^{{\rm skew}(0,j)}_{\C}$ defined below 
are denoted in \cite{JMT,JMMT} by 
$\F_{n,l}$, $\widehat{W}_{n,l}$,   
$W_{n,l}$ and 
$\widehat{\Zc}^{(j)}$, respectively.}. 
Let $\F^{\rm skew}_{\C,n,l}$
be the space of polynomials $P(X_1,\ldots,X_l)$ 
with coefficients in 
$\C[z_1^{\pm 1},\ldots,z_n^{\pm 1}]$, satisfying the following conditions:
\be
&&\mbox{$P$ is skew-symmetric in $X_1,\ldots,X_l$
(it is an empty condition when $l=0,1$),}\\
&&\deg_{X_i}P\le n-1.
\en
We denote by $\F^{\rm skew,\mathfrak S_n}_{\C,n,l}$ the
$\mathfrak S_n$-invariant subspace of $\F^{\rm skew}_{\C,n,l}$.
Let $W^{\rm skew}_{\C,n,l}$ denote the subspace of 
$\F^{\rm skew,\mathfrak S_n}_{\C,n,l}$ consisting of elements $P$ such that 
\be
&&P\vert_{X_1=z_1{}^{-1}=-z_2{}^{-1}}=0\quad \mbox{for $n\ge2$ and $l\ge1$.}
\en
We set 
\be
W^{{\rm skew}\ge0}_{\C,n,l}
:=W^{\rm skew}_{\C,n,l}\cap \C[z_1,\ldots,z_n][X_1,\ldots,X_l].
\en
Let further $\widehat{\Zc}^{{\rm skew}(i,j)}_{\C}[m]$ 
denote the space of sequences 
$(P_{n,l})_{n-2l=m}\in\prod_{\substack{n\ge 0\\n-2l=m}}
\F^{\rm skew}_{\C,n,l}$, satisfying the conditions
\be
&&P_{n+2,l+1}(X_1,\ldots,X_l,z^{-1}|z_1,\ldots,z_n,z,-z)\\
&&\quad =z^{-n-1+i}\prod_{a=1}^l(1-X_a^2z^2)\cdot
P_{n,l}(X_1,\ldots,X_l|z_1,\ldots,z_n).
\en
We set 
$\F^{\rm skew,\mathfrak S_n}_{\C,n}=\mathop\oplus\limits_{l=0}^n\F^{\rm skew,\mathfrak S_n}_{\C,n,l}$,
$W^{\rm skew}_{\C,n}=\mathop\oplus\limits_{l=0}^nW^{\rm skew}_{\C,n,l}$, 
$\widehat{\Zc}^{{\rm skew}(i,j)}_{\C}
=\mathop\oplus\limits_{\substack{m\in\Z\\m\equiv i-j\bmod 2}}
\widehat{\Zc}^{{\rm skew}(i,j)}_{\C}[m]$.

The spaces $\F^{\rm skew,\mathfrak S_n}_{\C,n,l}$, $W^{\rm skew}_{\C,n,l}$,
$\widehat{\Zc}^{{\rm skew}(i,j)}_{\C}$    
admit an action of $U_{\sqrt{-1}}$ (see \cite{JMT,JMMT}). 
As noted in \eqref{eq:Pnl},
we have an embedding of $U_{\sqrt{-1}}$-modules 
\bea\label{eq:Fi}
&&
\begin{array}{rcl}
\iota:
\bigl(\F^{\mathfrak S_n}_{A,n,l}\bigr)_{\sqrt{-1}}&\longrightarrow&
\F^{\rm skew,\mathfrak S_n}_{\C,n,l},\\[5pt]
P&\mapsto &c_{n,l} P\cdot
\prod_{j<j'}\frac{X_j-X_{j'}}{X_j+X_{j'}},
\end{array}
\ena
where $c_{n,l}\in\C\backslash\{0\}$. 
We can choose $c_{n,l}$ so that we have a map
\be
\iota:(\widehat{\Zc}^{(i,j)}_A)_{\sqrt{-1}}
\longrightarrow
\widehat{\Zc}^{{\rm skew}(i,j)}_{\C}, 
\en
and that when composed with the morphism 
\be
V_{\sqrt{-1}}(\Lambda_i)\otimes V_{\sqrt{-1}}(-\Lambda_j)
\longrightarrow (\widehat{\Zc}^{(i,j)}_A)_{\sqrt{-1}}
\en
the following are valid.
\be
&&v_0\otimes\vb_0\mapsto{\bf 1}^{(0,0)}_{\sqrt{-1}}=(1,X,X\wedge X^3,\ldots)
\in\widehat{\Zc}^{{\rm skew}(0,0)}_{\C}[0],\\
&&v_0\otimes\vb_1\mapsto{\bf 1}^{(0,1)}_{\sqrt{-1}}
=(1,X^2,X^2\wedge X^4,\ldots)
\in\widehat{\Zc}^{{\rm skew}(0,1)}_{\C}[1],\\
&&v_1\otimes\vb_0\mapsto{\bf 1}^{(1,0)}_{\sqrt{-1}}=(1,X,X\wedge X^3,\ldots)
\in\widehat{\Zc}^{{\rm skew}(1,0)}_{\C}[1],\\
&&v_1\otimes\vb_{-1}\mapsto{\bf 1}^{(1,1)}_{\sqrt{-1}}
=(1,1,1\wedge X^2,\ldots)
\in\widehat{\Zc}^{{\rm skew}(1,1)}_{\C}[0].
\en
Here we used the wedge product notation 
\be
&&P_1\wedge P_2:=\frac{1}{l_1!l_2!}
\Skew P_1(X_1,\ldots,X_{l_1})
P_2(X_{l_1+1},\ldots,X_{l_1+l_2}),\\
&&(\Skew f)(X_1,\ldots,X_l)=
\sum_{\sigma\in \mathfrak S_l}
(\sgn\sigma) f(X_{\sigma(1)},\ldots,X_{\sigma(l)}). 
\en

The following result was proved in \cite{JMT,JMMT}. 
\begin{prop}\label{prop:Wskew-is-cyclic}
We have
\be
W^{\rm skew}_{\C,n}=U_{\sqrt{-1}}1_n,
\quad 
W^{{\rm skew}\ge 0}_{\C,n}=U^{\ge 0}_{\sqrt{-1}}1_n.
\en
\end{prop}
Let us finish the proof of Proposition \ref{prop:W-is-cyclic}.

\medskip
\noindent{\it Proof of Proposition \ref{prop:W-is-cyclic}}.  
\quad
In order to show the
equality $U^{\ge 0}_{q}1_n=W^{\ge 0}_{n}$, it is enough to show the
equality of their characters. 
We have an inclusion of the $A$-modules
$U^{\ge 0}_A1_n\subset W^{\ge 0}_{A,n}$.  
They are free $A$-modules
because both of them are $A$-submodules of the free $A$-module
$A[z_1,\ldots,z_n][X_1,\ldots,X_l]$. 
Using Proposition
\ref{prop:Wskew-is-cyclic}, we have
\be
&&{\rm ch}\,W^{{\rm skew}\ge 0}_{\C,n}=
{\rm ch}\,U^{\ge 0}_{\sqrt{-1}}1_n\leq
{\rm ch}\,(U^{\ge 0}_A1_n)_{\sqrt{-1}}\\
&&={\rm ch}\,U^{\ge 0}_A1_n\leq
{\rm ch}\,W^{\ge 0}_{A,n}=
{\rm ch}\,(W^{\ge 0}_{A,n})_{\sqrt{-1}}\leq
{\rm ch}\,W^{{\rm skew}\ge 0}_{\C,n}.
\en
Here the first inequality follows from the surjective map
$(U^{\ge 0}_A1_n)_{\sqrt{-1}}\epi U^{\ge 0}_{\sqrt{-1}}1_n$,
and the last inequality follows from the injective map
$(W^{\ge 0}_{A,n})_{\sqrt{-1}}\mono W^{{\rm skew}\ge 0}_{\C,n}$
induced by \eqref{eq:Fi}.
In particular, we have
${\rm ch}\,U^{\ge 0}_A1_n={\rm ch}\,W^{\ge 0}_{A,n}$,
which implies
$U^{\ge 0}_{q}1_n=W^{\ge 0}_{n}$.
Noting that $(z_1\cdots z_n)^{-L}1_n\in U_{q}1_n$, 
we obtain
\be
W_{n}=\cup_{L}(z_1\cdots z_n)^{-L}W^{\ge 0}_{n}
=\cup_{L}(z_1\cdots z_n)^{-L}U^{\ge 0}_{q}1_n
=U_{q} 1_n.
\en
\qed

Arguing similarly as in the previous subsection, 
we obtain the following isomorphism
conjectured in \cite{JMMT}.  
\begin{thm}\label{thm:VVZ}
\bea
&&V_{\sqrt{-1}}(\Lambda_i)\ftens
V_{\sqrt{-1}}(-\Lambda_j)
\simeq\widehat{\Zc}^{{\rm skew}(i,j)}_{\C}
\ena
\end{thm}
\subsection{Characters}
In what follows we use the standard symbol 
$(v)_n:=\prod_{j=1}^n(1-v^j)$. 
The character of $W^{{\rm skew}\ge0}_{\C,n,l}$ 
was computed by Nakayashiki \cite{N}. 
{}From the results of the previous subsection, 
we conclude that $W^{\ge 0}_{n}$ has the same character:
\begin{cor}\label{CHAR+}
We have
\be
\ch_{v,z}W^{\ge 0}_n=\ch_{v,z}U_q^{\ge 0}v_+^{\otimes n}
=\sum_{l=0}^n\frac{z^{n-2l}}{(v)_l(v)_{n-l}}.
\en
\end{cor}

In particular, taking the sum over $n$ we obtain the 
known character formula of the integrable $\widehat{\mathfrak{sl}}_2$-modules
of level $-1$ \cite{Melzer}:
\bea\label{MEL}
\ch_{v,z}\bigl(V(-\Lambda_j)\bigr)
=\sum_{n\equiv j\bmod 2}v^{(n^2-j)/4}
\frac{\ch_{v,z}\pi_1^{*n}}{(v)_n},
\ena
where
\be
\ch_{v,z}\pi_1^{*n}=\sum_{l=0}^n
\frac{(v)_n}{(v)_l(v)_{n-l}}z^{n-2l}
\en
is the graded character of the fusion product of $n$ 
copies of $2$-dimensional irreducible $\mathfrak{sl}_2$-modules,
see, e.g., (2.11) in \cite{FJLM}.

In this paper
we considered the filtration of the tensor product of the level $1$ and $-1$ modules. 
The graded space associated with the induced filtration \eqref{induced-filtration} 
of $V(-\Lambda_{j})$ has the character \eqref{MEL}. 
It is known that similar filtrations of tensor products 
exist in a very general setting, see \cite{BN} and Section A.2.
In general, however, the subspaces defining the filtration
are not generated by tensor products of extremal vectors. 
In the case of 
integrable $\widehat{\mathfrak{sl}}_2$-modules of level $-k$, 
the following fermionic formula is known (the formula (2.14) in \cite{FJLM} in the limit $N \to \infty$):
\be
\ch_{v,z}V\bigl(-(k-j)\Lambda_0-j\Lambda_1\bigr)
=\sum_{n\equiv j\bmod 2}K^{(k)}_{j,(n,\underbrace{\scriptstyle{0,\ldots,0}}_{k-1})}(v)
\frac{ch_{v,z}\pi_1^{*n}}{(v)_n}.
\en
Here $K^{(k)}_{j,(n,0,\ldots,0)}(v)$ 
denotes the level restricted Kostka polynomial for 
$\mathfrak{sl}_2$, see e.g. (2.9) in \cite{FJLM}.
We conjecture that the right hand side gives the character of the 
associated graded space of the induced filtration mentioned above. 

\appendix

\section{Crystal and global bases}
\subsection{Summary of known facts}
In this subsection, we briefly summarize some of the basic 
definitions and results on crystal and global bases
which are relevant to the main text. 
We also give proofs of Propositions 
\ref{prop:intersection} and \ref{CORO}
at the end. 
Since our application in this paper is limited to the
case of $U_q(\widehat{\mathfrak{sl}}_2)$, 
we do not spare time to prepare
the notation for the general case. 
However, most of the statements are 
valid for arbitrary quantum affine 
algebras under suitable modifications.

Let $R$ be a subring of $K$. We use $R=A$, $A_0$ or $A_\infty$,
where $A_0$ (resp., $A_\infty$) is the ring of rational functions in $q$
which are regular at $q=0$ (resp., $q=\infty$).
Let $V$ be a vector space over $K$. An $R$-submodule $L\subset V$
is called an $R$-lattice of $V$ if $L$ is $R$-free 
and $V=K\otimes_RL$.

Let $V$ be an integrable $U_q$-module. In particular, $V$ has a weight space
decomposition, $V=\oplus_{\lambda\in P}V_\lambda$. 
The operators $\tilde f_i$ and
$\tilde e_i$ are defined as usual (see (2.2.2) in \cite{K1}).
A crystal base of a $U_q$-module $V$ is a pair $(L(V),B(V))$ of
an $A_0$-lattice $L(V)$ and a basis $B(V)$ of the $\C$-vector space $L(V)/qL(V)$
satisfying the following conditions:
\bea
&&\tilde e_iL(V)\subset L(V),\quad\tilde f_iL(V)\subset L(V)\hbox{ for any }i,\\
&&\tilde e_iB(V)\subset B(V)\sqcup\{0\}, 
\quad\tilde f_iB(V)\subset B(V)\sqcup\{0\},\\
&&L(V)=\oplus_{\lambda\in P}L(V)_\lambda,\quad
B(V)=\bigsqcup_{\lambda\in P}B(V)\cap(L(V)/qL(V))_\lambda,\\
&&b'=\tilde f_i b\hbox{ if and only if }b=\tilde e_ib'
\hbox{ for all }b,b'\in B(V)\hbox{ and }i.
\ena

Let $V_A$ be an $A$-lattice of $V$, $L_0$ an $A_0$-lattice, and $L_\infty$
an $A_\infty$-lattice. Set $E=L_0\cap L_\infty\cap V_A$.
The triplet $(L_0,L_\infty,V_A)$ 
is called {\em balanced} if the mapping
of $\C$-vector spaces
\bea
E\longrightarrow L_0/qL_0
\ena
is an isomorphism. We denote by $G$ the inverse map of this isomorphism.
Suppose that $(L_0,B)$ is a crystal base of $V$. 
The basis $\{G(b)\mid b\in B\}$ of $V$ 
is called a {\em global basis}. We have $V=\oplus_{b\in B}KG(b)$
and $V_A=\oplus_{b\in B}AG(b)$.

There is an involution of $U_q$ called the bar involution:
\bea
\overline q=q^{-1},\quad
\overline e_i=e_i,\quad
\overline f_i=f_i,\quad
\overline{q^h}=q^{-h}.
\ena
Let $V$ be a $U_q$-module. 
An involution $\overline{\phantom{u}}$ of $V$ is called
a bar involution if 
$\overline{a v}=\overline{a}\overline{v}$ holds for
$a\in U_q$, $v\in V$.  

An extremal module $V(\lambda)$ ($\lambda\in P$) admits 
a bar involution such that $\overline{u_\lambda}=u_\lambda$. 
We take $V_A(\lambda)=U_Au_\lambda$ as its $A$-lattice. 
There exists a crystal base
$(L(\lambda),B(\lambda))$ of $V(\lambda)$ such that the triple
$(L(\lambda),\overline L(\lambda),V_A(\lambda))$ is balanced.
The construction is as follows.

First, suppose that $\lambda\in P_+$. We define
\bea
&&L(\lambda)=\sum_{m=0}^\infty\sum_{i_1,\ldots,i_m\in I}
A_0\tilde f_{i_1}\cdots\tilde f_{i_m}u_\lambda,\\
&&B(\lambda)=
\{\tilde f_{i_1}\cdots\tilde f_{i_m}u_\lambda\in L(\lambda)/qL(\lambda)\mid
m\in\Z_{\ge 0},i_1,\ldots,i_m\in I\}\backslash\{0\}.
\ena
\begin{prop}\cite{K0}
For $\lambda\in P_+$, the pair $(L(\lambda),B(\lambda))$ is a crystal base
of $V(\lambda)$. 
The triplet $(L(\lambda),\overline{L(\lambda)},V_A(\lambda))$
is balanced.
\end{prop}
\noindent
Similarly, we can construct the crystal and global bases for $V(\lambda)$ when
$\lambda\in P_-$. 
By abuse of notation we use $u_\lambda\in B(\lambda)$.

If $\lambda,\mu\in P_+$, we have an embedding
\bea
V(\lambda+\mu)\simeq U_q(u_\lambda\otimes u_\mu)
\subset V(\lambda)\otimes V(\mu)
\ena
such that $u_{\lambda+\mu}\mapsto u_\lambda\otimes u_\mu$.
It induces an embedding of crystal
$B(\lambda+\mu)\subset B(\lambda)\otimes B(\mu)$.

Now, consider the tensor product $V(\lambda)\otimes V(-\mu)$
where $\lambda,\mu\in P_+$. The vector $u_\lambda\otimes u_{-\mu}$
is a cyclic vector of $V(\lambda)\otimes V(-\mu)$. In fact, we have
\bea
U_A(u_\lambda\otimes u_{-\mu})=V_A(\lambda)\otimes V_A(-\mu).
\ena
There exists a unique bar involution of 
$V(\lambda)\otimes V(-\mu)$ such that 
$\overline{u_\lambda\otimes u_{-\mu}}=u_\lambda\otimes u_{-\mu}$. 
In general, $\overline{u\otimes v}$ is not 
equal to $\overline u\otimes\overline v$. 
However, we have
\bea
\overline{u_\lambda\otimes v}=u_\lambda\otimes\overline v,\quad
\overline{v\otimes u_{-\mu}}=\overline v\otimes u_{-\mu}.
\ena
\begin{prop}\cite{L}
The pair $(L(\lambda)\otimes L(-\mu),B(\lambda)\otimes B(-\mu))$
is a crystal base of $V(\lambda)\otimes V(-\mu)$, 
and the triplet
$(L(\lambda)\otimes L(-\mu),
\overline{L(\lambda)\otimes L(-\mu)},
V_A(\lambda)\otimes V_A(-\mu))$ is balanced.
\end{prop}

In general, it is not true that $G(b_1\otimes b_2)=G(b_1)\otimes G(b_2)$
for $b_1\in B(\lambda)$ and $b_2\in B(-\mu)$. However, we have
\bea
&&G(u_\lambda\otimes b)=u_\lambda\otimes G(b)\hbox{ for any }b\in B(-\mu),\\
&&G(b\otimes u_{-\mu})=G(b)\otimes u_{-\mu}\hbox{ for any }b\in B(\lambda).
\ena

Let $\lambda\in P$, and write it as $\lambda=\xi-\eta$ where $\xi,\eta\in P_+$.
We have a surjection
\bea\label{XIETA}
p_{\xi,\eta}:V(\xi)\otimes V(-\eta)\rightarrow V(\lambda)
\ena
sending $u_\xi\otimes u_{-\eta}$ to $u_\lambda$. 
We set $L(\lambda)=p_{\xi,\eta}(L(\xi)\otimes L(-\eta))$. 
The map
$p_{\xi,\eta}$ induces
$\overline p_{\xi,\eta}:(L(\xi)/qL(\xi))
\otimes(L(-\eta)/qL(-\eta))
\rightarrow L(\lambda)/qL(\lambda)$. 
We set
$B(\lambda)=\overline p_{\xi,\eta}(B(\xi)\otimes B(-\eta))\backslash\{0\}$.
\begin{prop}\cite{L,K1}\label{HL}
The pair $(L(\lambda),B(\lambda))$ is a crystal base of $V(\lambda)$,
and the triplet $(L(\lambda),\overline{L(\lambda)},V_A(\lambda))$
is balanced. For $b\in B(\xi)\otimes B(-\eta)$ we have
\bea
&&p_{\xi,\eta}(G(b))=
G(\overline p_{\xi,\eta}(b)).
\ena
\end{prop}

The crystal base $(L(\lambda),B(\lambda))$ and the global base $G(b)$
$(b\in B(\lambda))$ are independent of the choice of $\xi,\eta\in P_+$
such that $\lambda=\xi-\eta$. In fact, they are obtained from a universal
object called the modified quantized enveloping algebra.
The modified quantized enveloping algebra $\tilde U_q$ is
\bea
\tilde U_q\seteq \oplus_{\lambda\in P}U_qa_\lambda\hbox{ where }
U_qa_\lambda=U_q/\sum_{h\in P^*}U_q(q^h-q^{\langle h,\lambda\rangle}).
\ena
For any $\xi,\eta\in P_+$ such that $\lambda=\xi-\eta$, 
we denote by
$\Phi_{\xi,\eta}:U_qa_\lambda\rightarrow V(\xi)\otimes V(-\eta)$ 
the $U_q$-linear mapping which sends $a_\lambda$ to
$u_\xi\otimes u_{-\eta}$.
\begin{prop}
There exists a unique $A_0$-lattice $\tilde L_\lambda$ of 
$U_qa_\lambda$ and a unique basis $\tilde B_\lambda$ of
$\tilde L_\lambda/q\tilde L_\lambda$ $(\lambda\in P)$
satisfying the following properties.

\bi
\item
The triplet $(\tilde L_\lambda,\overline{\tilde L_\lambda},U_Aa_\lambda)$
is balanced.

\item
The image of $\tilde L_\lambda$ by $\Phi_{\xi,\eta}$ is equal to
$L(\xi)\otimes L(-\eta)$.

\item
Let $\overline\Phi_{\xi,\eta}:\tilde L_\lambda/q\tilde L_\lambda
\rightarrow L(\xi)\otimes L(-\eta)/qL(\xi)\otimes L(-\eta)$
be the induced map. 
Then $\overline{\Phi}_{\xi,\eta}$ gives a bijection between 
$\{b\in\tilde B_\lambda\mid
\overline\Phi_{\xi,\eta}(b)\not=0\}$ and 
$B(\xi)\otimes B(-\eta)$. 
For $b\in\tilde B_\lambda$ we have
\bea
\Phi_{\xi,\eta}(G(b))=
G(\overline\Phi_{\xi,\eta}(b)).
\ena

\item
The set $\tilde B_\lambda$ has a structure of crystal such that
$B(\xi)\otimes B(-\eta)\sqcup\{0\}\subset\tilde B_\lambda\sqcup\{0\}$
is an embedding which commutes with the action of
$\tilde e_i,\tilde f_i$ $(i\in I)$.

\item
The set $\tilde B_\lambda$ is equal to the inductive limit
$\mathop{\varinjlim}\limits_{\xi,\eta\rightarrow\infty}
B(\xi)\otimes B(-\eta)\sqcup\{0\}$,
where we use the dominance ordering in $P_+$ in taking the limit.
\ei
\end{prop}
We write $\widetilde{B}_\lambda$ for $B(U_qa_\lambda)$.

\begin{prop}\label{GBTENSOR}
Let $\Phi_\lambda\colon U_qa_\lambda\to V(\lambda)$ be the surjective 
morphism sending $a_{\lambda}$ to $u_{\lambda}$. 
Then the induced morphism
$\overline\Phi_\lambda\colon \tilde L_\lambda/q\tilde L_\lambda
\to L(\lambda)/qL(\lambda)$ satisfies
$\overline\Phi_\lambda(\tilde B_\lambda)\subset B(\lambda)\sqcup\{0\}$.
Moreover
$\{b\in\tilde B_\lambda\,;\,\overline\Phi_\lambda(b)\not=0\}\to B(\lambda)$
is bijective, and
$\Phi_\lambda(G(b))=G(\overline\Phi_\lambda(b))$.
\end{prop}
Let $\mu\in P$. Suppose that $-w\mu\in P_+$ for some $w\in W$.
Namely, the weight $\mu$ is an extremal weight in the weight space of
the lowest weight module $V(w\mu)$. We identify, $V(\mu)\simeq V(w\mu)$,
$B(\mu)\simeq B(w\mu)$, and the extremal vector
$S_{w^{-1}}u_{w\mu}\in V(w\mu)$ with $u_\mu\in V(\mu)$.
\begin{prop}
There exists a subset $B^+_\mu$ of $B(\mu)$ such that
\bea
U^+_q u_\mu=\oplus_{b\in B^+_\mu}KG(b).
\ena
\end{prop}

\begin{prop}\label{APPPRO}
Suppose that $\mu=\xi-\eta$ where $\xi,\eta\in P_+$.
There exists an embedding of crystals
\bea
B(\mu)\subset B(\xi)\otimes B(-\eta)
\ena
such that
\bea\label{SUCHTHAT}
B^+_\mu=u_\xi\otimes B(-\eta).
\ena
\end{prop}

Let $\lambda\in P_+$ and $\mu\in P$.
By Proposition \ref{HL}, the tensor product $V(\lambda)\otimes V(\mu)$
has the global base $G(b)$ where $b\in B(\lambda)\otimes B(\mu)$.
Consider the submodule
$N_{\lambda,\mu}=U_q(u_\lambda\otimes u_\mu)\subset V(\lambda)\otimes V(\mu)$.
\begin{prop}\label{NONCYC}
There exists a subset $B_{\lambda,\mu}$ of $B(\lambda)\otimes B(\mu)$ such that
\bea
N_{\lambda,\mu}=\oplus_{b\in B_{\lambda,\mu}}KG(b).
\ena
The set $B_{\lambda,\mu}\sqcup\{0\}$ is invariant by $\tilde f_i,\tilde e_i$.
\end{prop}
This proposition follows from
\begin{lem}\cite{KL}\label{APPLEM}
Let $M$ be an integrable module with a crystal base $(L(M),B(M))$
and a bar involution. Let $M_A$ be an $A$-lattice of $M$ such that
$(L(M),\overline{L(M)},M_A)$ 
is balanced. Let $N^+$ be
a $U^+$-submodule of $M$ such that
\bea
N^+=\oplus_{b\in B_{N^+}}KG(b)\hbox{ for a subset }
B_{N^+}\subset B(M).
\ena
Set
\bea
B_N=\{\tilde f_{i_1}\cdots\tilde f_{i_m}b\mid
b\in B_{N^+}\}\backslash\{0\}.
\ena
Then, we have
\bea
U_qN^+=\oplus_{b\in B_N}KG(b).
\ena
\end{lem}

Here we give proofs of Propositions \ref{prop:intersection}
and \ref{CORO}.
\medskip

\noindent
{\it Proof of Proposition \ref{prop:intersection}.}
We take $\lambda=\Lambda_i$ and $\mu=\mu_n=\wt\vb_{-(n-i)}$
where $n\equiv i-j\bmod2$ in Proposition \ref{NONCYC}.
Let us prove that
\bea\label{LETUS}
B_{\lambda,\mu}\cap(u_\lambda\otimes B(\mu))=u_\lambda\otimes B^+_\mu.
\ena

We take $\xi,\eta$ as in Proposition \ref{APPPRO}. Then, we have
$B(\mu)\subset B(\xi)\otimes B(-\eta)$.
We can choose $\xi,\eta$ in such a way that
if $\langle h_i,\lambda\rangle=0$ then $\langle h_i,\xi\rangle=0$.
Note that if $\langle h_i,\xi\rangle=0$ we have $\tilde f_i u_\xi=0$.
Therefore, from the tensor product rule for $\tilde f_i$
(see \cite{K0}, (2.4.3)), we have
\bea\label{FACTION}
\tilde f_i(u_\lambda\otimes u_\xi)=
\begin{cases}
\tilde f_iu_\lambda\otimes u_\xi&\hbox{if $\langle h_i,\lambda\rangle>0$};\\
0&\hbox{otherwise}.
\end{cases}
\ena

By Lemma \ref{APPLEM}, we know that $B_{\lambda,\mu}$ is obtained from
$u_\lambda\otimes B^+_\mu$ by applying $\tilde f_i$'s.
We have
\bea\label{4.21}
B^+_\mu=u_\xi\otimes B(-\eta)
\ena
by (\ref{SUCHTHAT}). Since
$u_\lambda\otimes u_\xi\in B(\lambda+\xi)$ and
$B(\lambda+\xi)\subset B(\lambda)\otimes B(\mu)$ is invariant by $\tilde f_i$,
we have $B_{\lambda,\mu}\subset B(\lambda+\xi)\otimes B(-\eta)$.

We have
\bea
B_{\lambda,\mu}\cap(u_\lambda\otimes B(\mu))\subset
(B(\lambda+\xi)\cap u_\lambda\otimes B(\xi))\otimes B(-\eta).
\ena
{}From (\ref{FACTION}) it follows that
$B(\lambda+\xi)\cap u_\lambda\otimes B(\xi)=u_\lambda\otimes u_\xi$.
Using (\ref{4.21}) we have (\ref{LETUS}).

Let us prove (\ref{prop:intersection}). It suffices to show that
\be
\bigcap_nB_{\lambda,\mu_n}=\{0\}.
\en

Suppose that $b_1\otimes b_2\in B_{\lambda,\mu_n}$. Let us show
that the actions of $\tilde e_i$'s bring $b_1\otimes b_2$ to
$u_\lambda\otimes b_2'\in B_{\lambda,\mu_n}$. If $b_1\not= u_\lambda$,
we have an 
$i$ such that $\tilde e_ib_1\not=0$. By the tensor product rule,
we see that there exists an $l>0$ such that
\be
(\tilde e_i)^l(b_1\otimes b_2)=\tilde e_ib_1\otimes(\tilde e_i)^{l-1}b_2\not=0.
\en
Since $B(\lambda)$ is connected, the assertion follows from this.

Now, we have $u_\lambda\otimes b_2\in B_{\lambda,\mu_n}$ for all $n$.
By (\ref{LETUS}) we have $b_2\in B_{\mu_n}$ for all $n$. Since
$\cap_nB_{\mu_n}=\{0\}$, we have the assertion (\ref{prop:intersection}).
\qed
\medskip

\noindent
{\it Proof of Proposition \ref{CORO}.}
We prove the first isomorphism. The second isomorphism then follows.
It suffices to prove
\be
U_A(v_i\otimes\vb_{n-i})\cap U_q(v_i\otimes\vb_{n+2-i})=
U_A(v_i\otimes\vb_{n+2-i}).
\en
Using the notation in the proof of Proposition \ref{prop:intersection},
we have 
\be
&&U_A(v_i\otimes\vb_{n-i})=\oplus_{b\in B_{\lambda,\mu_n}}AG(b),\\
&&U_q(v_i\otimes\vb_{n+2-i})=\oplus_{b\in B_{\lambda,\mu_{n+2}}}KG(b),\\
&&U_A(v_i\otimes\vb_{n+2-i})=\oplus_{b\in B_{\lambda,\mu_{n+2}}}AG(b).
\en
The assertion is clear from these equalities.
\qed
\medskip

\subsection{Filtration on $V(\xi)\otimes V(-\eta)$}

In the previous subsection, 
we prepared basic definitions
and results which are used in the main text in the construction of
the filtration of the tensor product of level $1$ and level $-1$
$U_q(\slth)$-modules. The submodules which constitute the filtration
in this case are generated by single vectors of the form
$v_i\otimes \vb_{n-i}$, and these vectors $v_n,\vb_n$ are the extremal
vectors. In a more general situation, i.e., $V(\xi)\otimes V(\eta)$
where $\xi,\eta\in P$ for affine quantum algebras other than $U_q(\slth)$
and/or level $k$ of $\xi,-\eta$ is greater than 1, the existence of
similar filtrations was proved in \cite{BN}. In the below,
we briefly state the construction.

Let $B(\infty)$ be the crystal of $U^-_q$ \cite{K0}.
It has a unique element $u_\infty$
which has weight zero, and is given in the form
\be
B(\infty)=\{\tilde f_{i_1}\cdots\tilde f_{i_n}u_\infty\}\backslash\{0\}.
\en
For $b\in B(\infty)$ we have
$\varepsilon_i(b)={\rm max}\{n\mid \tilde e_i^nb\not=0\}\geq0$ and
$\langle h_i,\wt b\rangle+\varepsilon_i(b)=\varphi_i(b)$.
Note that $\varphi_i(b)$ is finite and can be negative.
Similarly, we have the crystal $B(-\infty)$ of $U^+_q$:
\be
B(-\infty)=\{\tilde e_{i_1}\cdots\tilde e_{i_n}u_{-\infty}\}\backslash\{0\}.
\en
The weight of $u_{-\infty}$ is zero and
$\varphi_i(b)={\rm max}\{n\mid \tilde f_i^nb\not=0\}\geq0$.

Let $\lambda\in P$. We denote by $T_\lambda$ the crystal consisting of
a single element $t_\lambda$ such that
$\varepsilon_i(t_\lambda)=\varphi_i(t_\lambda)=-\infty$. Suppose that
$b$ is an element of a crystal such that $\varphi_i(b)$ is finite.
The equality $\varepsilon_i(t_\lambda)=-\infty$ implies
that $x(b\otimes t_\lambda)=xb\otimes t_\lambda$ for $x=\tilde f_i$ or
$\tilde e_i$. Similarly, if $\varepsilon_i(b)$ is finite
we have $x(t_\lambda\otimes b)=t_\lambda\otimes xb$.

We have an isomorphism of crystals
\be
B(U_qa_\lambda)\simeq B(\infty)\otimes T_\lambda\otimes B(-\infty).
\en
The right hand side has a decomposition into the crystals $B(\zeta)$ ($\zeta\in P$).
The decomposition is appropriately described by introducing
the star crystal structure on
$B(\tilde U_q)=\sqcup_{\lambda\in P}B(U_qa_\lambda)$.

Let $*$ be an anti-involution of $U_q$ given by
\be
q^*=q,\quad e_i^*=e_i,\quad f_i^*=f_i,\quad(q^h)^*=q^{-h}.
\en
We define $*$ on $\tilde U_q$ by $a_\lambda^*=a_{-\lambda}$.
It induces involutions of $B(\infty)$, $B(-\infty)$ and  
$B(\tilde U_q)$,
which we also denote by $*$. We have
\be
(b_1\otimes t_\lambda\otimes b_2)^*=
b_1^*\otimes t_{-\lambda-\wt b_1-\wt b_2}\otimes b_2^*.
\en
We define crystal structures on $B(\infty)$, $B(-\infty)$ and  
$B(\tilde U_q)$ by
\be
\tilde e_i^*=*\circ\tilde e_i\circ*,\quad
\tilde f_i^*=*\circ\tilde f_i\circ*.
\en
We call the actions of $\tilde e_i^*$ and $\tilde f_i^*$
the star crystal actions. They commute with the original crystal actions
$\tilde e_j$ and $\tilde f_j$. We have the structure of 
bi-crystal on
$B(\infty)$, $B(-\infty)$ and $B(\tilde U_q)$. 
We have $\wt^*(b)=\wt(b^*)$.
In particular, $\wt^*(U_qa_\lambda)=-\lambda$.
or equivalently, $\wt^*(b_1\otimes t_\lambda\otimes b_2)=-\lambda$.
See \cite{K2} for other useful formulas
of the star crystal actions on $B(\tilde U_q)$.

Let $\mathfrak g=\mathfrak{sl}_2$, and denote the simple root by
$\alpha=2\varpi$. It is easy to see that we have decompositions
\be
&&B(\infty)\otimes T_{n\varpi}\otimes B(-\infty)\\
&&\hspace{70pt}=\begin{cases}
B(n\varpi)\sqcup B((n+2)\varpi)\sqcup B((n+4)\varpi)\sqcup\dots
&\hbox{if $n\geq0$;}\\
B(n\varpi)\sqcup B((n-2)\varpi)\sqcup B((n-4)\varpi)\sqcup\dots
&\hbox{if $n<0$,}
\end{cases}
\en
as crystals, forgetting the star crystal actions. For example,
$\{u_\infty\otimes t_0\otimes u_{-\infty}\}$ is isomorphic to
$B(0)$ in the crystal action. In fact, it is also isomorphic to
$B(0)$ in the star crystal action. Therefore, the subset
$\{u_\infty\otimes t_0\otimes u_{-\infty}\}$ of
$B(\infty)\otimes T_0\otimes B(-\infty)$ is isomorphic to $B(0)\times B(0)$
as bi-crystal.
For bi-crystals, we use the symbol $\times$. It is not a product of crystals.
The crystal structure of the second component represents the star crystal
structure.

Similarly, $B(\varpi)\times B(-\varpi)$ is contained
in the union $B(\infty)\otimes T_\varpi\otimes B(-\infty)\cup
B(\infty)\otimes T_{-\varpi}\otimes B(-\infty)$. The identification is
\be
(u_\varpi,u_{-\varpi})&=&u_\infty\otimes t_\varpi\otimes u_{-\varpi},\\
(u_{-\varpi},u_{-\varpi})
&=&\tilde fu_\infty\otimes t_\varpi\otimes u_{-\varpi},\\
(u_{-\varpi},u_\varpi)&=&u_\infty\otimes t_{-\varpi}\otimes u_{-\varpi},\\
(u_\varpi,u_\varpi)&=&u_\infty\otimes t_\varpi\otimes \tilde eu_{-\varpi}.
\en

In general, we have the decomposition in the $\mathfrak{sl}_2$ case.
\be
B(\tilde U_q(\mathfrak{sl}_2))
=\sqcup_{n=0}^\infty B(n\varpi)\times B(-n\varpi),
\en
where we have the identification
$u_\infty\otimes t_{n\varpi}\otimes u_{-\infty}=(u_{n\varpi},u_{-n\varpi})$.

The decomposition of the bi-crystal $B(\tilde U_q(\mathfrak g))$
in the general case of quantum affine  algebras
was conjectured in \cite{K2}, 
and proved in \cite{BN}.
We denote by $B_0(\lambda)$ the connected component of $B(\lambda)$
which contains $u_\lambda$. There is an action of $W$ on
$\bigsqcup_{\lambda\in P}B(\lambda)\times B_0(-\lambda)$ 
induced by the isomorphism (\ref{ISOMOR}).
\begin{prop} We have
\bea
B(\tilde U_q)=\sqcup_{\lambda\in P}(B(\lambda)\times B_0(-\lambda))/W.
\ena
\end{prop}
In particular, we have
$u_\infty\otimes t_\lambda\otimes u_{-\infty}
=(u_\lambda,u_{-\lambda})$.

Now, we consider $\xi=\sum_i\xi_i\Lambda_i$ and
$\eta=\sum_i\eta_i\Lambda_i\in P_+$ such that
$\langle c,\xi\rangle=\langle c,\eta\rangle>0$. We set $\zeta=\xi-\eta$.
The level of $\zeta$ is zero. The crystal
$B(\xi)\otimes B(-\eta)$ is regarded as a subset of
\bea
B(U_qa_\zeta)=\sqcup_{\lambda\in P}
(B(\lambda)\otimes B_0(-\lambda)_{-\zeta})/W.
\ena
We have the characterization of $B(\xi)\otimes B(-\eta)$
by using the star crystal structure.
\begin{prop}
We have
\bea
&&B(\xi)\otimes B(-\eta)
=\{b\in B(U_qa_\zeta)\mid\varepsilon_i^*(b)\leq\xi_i,\varphi_i^*(b)\leq\eta_i
\hbox{ for all }i\in I\}.
\ena
\end{prop}
\begin{proof}
By Proposition \ref{Lu}, we have
\be
V(\xi)\otimes V(-\eta)\simeq U_qa_\zeta/\left(\sum_iU_qf^{1+\xi_i}a_\zeta
+\sum_iU_qe^{1+\eta_i}a_\zeta\right).
\en
Set $I_{\xi,\eta}=\sum_iU_qf^{1+\xi_i}+\sum_iU_qe^{1+\eta_i}$.
We have
\be
I_{\xi,\eta}a_\zeta=U_qa_\zeta\cap
\left(\sum_i\tilde U_qf^{1+\xi_i}+\sum_i\tilde U_qe^{1+\eta_i}\right).
\en
Therefore, we have
\be
(I_{\xi,\eta}a_\zeta)^*=a_{-\zeta}U_q
\cap\left(\sum_if_i^{1+\xi_i}\tilde U_q+\sum_ie_i^{1+\eta_i}\tilde U_q\right).
\en
{}From this follows
\be
I_{\xi,\eta}a_\zeta=
\oplus_{b\in B_{\xi,\eta}}KG(b),
\en
where
\be
B_{\xi,\eta}=\left\{b\in B(\tilde U_q)\,;\ \parbox{270pt}{there exists an 
$i\in I$ such that
$\varepsilon_i(b^*)\geq1+\xi_i$ or
$\varphi_i(b^*)\geq1+\eta_i$ holds}\right\}
\en
Since $\varepsilon_i^*(b)=\varepsilon_i(b^*)$
and $\varphi_i^*(b)=\varphi_i(b^*)$, the assertion follows.
\end{proof}
Since $\langle h_i,\wt^*(b)\rangle+\varepsilon_i^*(b)=\varphi_i^*(b)$,
the condition $\varepsilon_i(b^*)\leq\xi_i$ is equivalent to
$\varphi_i(b^*)\leq\eta_i$.

Let $\lambda,\lambda'\in P$ be of level zero, we consider a partial ordering
$\lambda\geq\lambda'$ if and only if
${\rm cl}(\lambda-\lambda')\in\sum_{i\in I\backslash\{0\}}\Z_{\geq0}\alpha_i$.
Here ${\rm cl}(\lambda)\in P/\Z\delta$ is the classical part of $\lambda$.
We fix a total ordering for $\lambda,\lambda'\in P$ which is a refinement of
the partial ordering. We denote it also by $\geq$.

We fix a set of representatives
$P_+^{(0)}\subset\{\lambda\in P\mid\langle c,\lambda\rangle=0\}$
with respect to the action of the Weyl group $W$ in such a way that
for any $\lambda\in P^+_0$, we have
${\rm cl}(\lambda)\in\sum_{i\in I\backslash\{0\}}{\Q_{\geq0}}\alpha_i$.
We denote the isotropy subgroup of $\lambda$ by $W_\lambda$.
Define a filtration
$F^\lambda_{\xi,\eta}$ $(\lambda\in P_+^{(0)})$ of
$V(\xi)\otimes V(-\eta)$ by
\bea
F^{\geq\lambda}_{\xi,\eta}=
\bigoplus_{\lambda'\geq\lambda,\lambda'\in P_+^{(0)}}
\left(\bigoplus_{b\in(B(\lambda')\times B_0(-\lambda')/W_{\lambda'})
\cap(B(\xi)\otimes B(-\eta))}KG(b)\right)
\ena
Similarly, we define $F^{>\lambda}_{\xi,\eta}$.

The following proposition follows from \cite{BN}.
\begin{prop}
The subspace $F^{\geq\lambda}_{\xi,\eta}$ is $U_q$-invariant.
The subquotient $F^{\geq\lambda}_{\xi,\eta}/F^{>\lambda}_{\xi,\eta}$
is isomorphic to a direct sum of copies of $V(\lambda)$.
\end{prop}

For $\mathfrak g=\widehat{\mathfrak{sl}}_2$, we have an explicit description of
$B_0(-\lambda)_{-\zeta}$ by using ``paths''. 
See \cite{Nakashima}.
Note that the isotropy subgroup $W_\lambda$ 
is trivial in this case.
It is enough to consider the following two cases.

{\it Case} 1: $\lambda=0$ and $\zeta=0$.

{\it Case} 2: $\lambda=n(\Lambda_1-\Lambda_0)+m\delta$ $(n>0)$,
$|\xi_1-\eta_1|\leq n$ and $\xi_1-\eta_1\equiv n\bmod2$.

\noindent
{\it Case} 1 is trivial.
In {\it Case} 2, we embed
\be
B_0(-\lambda)\simeq B_0(n(\Lambda_1-\Lambda_0))\otimes T_{-m\delta}
\subset (B(\Lambda_1-\Lambda_0)^{\otimes n})\otimes T_{-m\delta}.
\en
We have the identification
\be
B(\Lambda_1-\Lambda_0)=\{z^\mu v_\varepsilon\mid
\mu\in\Z,\varepsilon=\pm1\}.
\en
In this identification, an element
$b=z^{\mu_1}v_{\varepsilon_1}\otimes\dots\otimes z^{\mu_n}v_{\varepsilon_n}
\otimes t_{-m\delta}$ belongs to $B_0(-\lambda)_{-\zeta}$
if and only if
\bea
&&\mu_{l+1}-\mu_l=
\begin{cases}
0&\hbox{ if }(\varepsilon_l,\varepsilon_{l+1})=(+,+),(-,+),(-,-);\\
1&\hbox{ if }(\varepsilon_l,\varepsilon_{l+1})=(+,-).
\end{cases}\\
&&\mu_1+\dots+\mu_n=m,\\
&&\varepsilon_1+\dots+\varepsilon_n=\eta_1-\xi_1.
\ena
The conditions $\varepsilon^*_i(b)\leq\xi_i$ $(i=0,1)$ are equivalent to
the ``level restriction''
\bea
&&{\rm max}\{\varepsilon_1+\dots+\varepsilon_l;1\leq l\leq n\}\leq\xi_0,\\
&&{\rm min}\{\varepsilon_1+\dots+\varepsilon_l;1\leq l\leq n\}\geq\xi_1.
\ena

{\it Acknowledgments.}\quad
MJ is partially supported by 
the Grant-in-Aid for Scientific Research (B2) no.12440039, 
Japan Society for the Promotion of Science (JSPS),
MK is partially supported by 
Grant-in-Aid for Scientific Research (B1)13440006,
JSPS,
and TM is partially supported by 
(A1) no.13304010, JSPS,
EM is partially supported by the National Science
Foundation (NSF) grant DMS-0140460.
YT is supported by JSPS.


\end{document}